\documentclass[10pt,leqno,twoside,draft]{article}
\usepackage{amsmath}
\usepackage{amsfonts}
\usepackage{amsthm}
\usepackage{indentfirst}

\theoremstyle{plain}
\newtheorem{theorem}{\indent\rm T\,h\,e\,o\,r\,e\,m\;}[section]
\newtheorem{lemma}{\indent\rm L\,e\,m\,m\,a\;}[section]
\newtheorem{proposition}{\indent\rm P\,r\,o\,p\,o\,s\,i\,t\,i\,o\,n\;}[section]
\newtheorem{corollary}{\indent\rm C\,o\,r\,o\,l\,l\,a\,r\,y\;}[section]

\theoremstyle{definition}
\newtheorem{definition}{\indent\rm D\,e\,f\,i\,n\,i\,t\,i\,o\,n\;}[section]

\theoremstyle{remark}
\newtheorem{remark}{\indent\rm R\,e\,m\,a\,r\,k\;}[section]

\theoremstyle{example}
\newtheorem{example}{\indent\rm E\,x\,a\,m\,p\,l\,e\;}[section]

\usepackage[english]{babel}
\usepackage{pdfsync}
\usepackage{graphicx}
\usepackage{pgfplots}

\newcommand\vv{\textsc{v}}

\newcommand{\GG}{\mathcal G}
\newcommand{\G}{\mathcal G}
\newcommand{\erre}{\mathbb{R}}
\newcommand{\Hmu}{{H_\mu^1}}
\newcommand{\R}{{\mathbb R}}
\newcommand{\f}{\frac}
\newcommand{\E}{\mathcal E}
\newcommand{\EE}{\mathcal E}

\tikzstyle{nodino}=[circle,draw,fill,inner sep=0pt,minimum size=0.5mm]
\tikzstyle{nodo}=[circle,draw,fill,inner sep=0pt,minimum size=1.5mm]
\tikzstyle{infinito}=[circle,inner sep=0pt,minimum size=0mm]

\title{Nonlinear dynamics on branched structures and networks}

\author{Riccardo Adami\thanks{Author partially supported by the FIRB 2012 project ``Dispersive dynamics: Fourier
Analysis and Variational Methods".}, Enrico Serra, Paolo Tilli \\ \ \\{\small  Dipartimento di Scienze
Matematiche ``G.L. Lagrange'', Politecnico di Torino } \\ {\small
Corso Duca degli Abruzzi, 24, 10129 Torino, Italy}}

\begin{document}
\thispagestyle{empty}


\begin{center}
{\sc\large Riccardo Adami,} \ {\sc\large Enrico Serra}
\ {\small and}
 \ {\sc\large Paolo Tilli}
\end{center}
\vspace {1.5cm}

\centerline{\large{\textbf{Nonlinear dynamics on branched structures and networks}}}
\renewcommand{\thefootnote}{\fnsymbol{footnote}}


                                                               %
\makeatletter                                                  %
\renewcommand*{\@seccntformat}[1]{
  \csname the#1\endcsname\;-                                   %
}                                                              %
\renewcommand{\section}{\@startsection{section}{1}{0mm}        %
   {1.5\baselineskip}
   {1\baselineskip}
   {\indent\normalfont\normalsize\bfseries}
   }                                                           %
\renewcommand*{\@seccntformat}[1]{
  \normalfont\bfseries\csname the#1\endcsname\;-               %
}                                                              %
\renewcommand\subsection{\@startsection                        %
  {subsection}{2}{0mm}
  {1.5\baselineskip}
  {1\baselineskip}
  {\indent\normalfont\normalsize\itshape}}
\renewcommand*{\@seccntformat}[1]{
  \normalfont\bfseries\csname the#1\endcsname\;-               %
}                                                              %
\renewcommand\subsubsection{\@startsection                     %
  {subsubsection}{2}{0mm}
  {1.5\baselineskip}
  {1\baselineskip}
  {\indent\normalfont\normalsize\texttt}}
\makeatother                                                   %
                                                               %

\renewcommand{\thefootnote}{\arabic{footnote}}
\setcounter{footnote}{0}

\vspace{1,5cm}
\begin{center}
\begin{minipage}[t]{10cm}

\small{ \noindent \textbf{Abstract.}
In these lectures we review on a recently developed line of research, concerning
the existence of ground states with prescribed mass (i.e. $L^2$-norm) for the {\em focusing} nonlinear Schr\"odinger
equation with a power  nonlinearity, {\em on noncompact quantum graphs.}

Nonlinear dynamics on graphs  has rapidly become a topical issue with many physical applications, ranging from
nonlinear optics to Bose-Einstein condensation. Whenever in a
physical experiment a ramified structure is involved (e.g. in the propagation of signals, in a circuit of quantum wires or
in trapping a boson gas), it can prove useful to approximate such a structure by a {\em metric graph}, or {\em network}.

For the Schr\"odinger equation 
it turns out that the 
{\em sixth} power in the nonlinear term of the energy (corresponding to the {\em quintic} nonlinearity in the evolution equation) is {\em critical}
 in the
sense that below that power the constrained energy is lower bounded irrespectively of the value of the mass ({\em subcritical case)}. On the other hand, if the nonlinearity
power equals six, then the lower boundedness depends on the value of the mass: below a {\em critical mass}, the 
constrained energy is lower bounded, beyond it, it is not.
For powers larger than six the 
constrained energy functional is never lower bounded, so that it is meaningless to speak about ground states ({\em supercritical case}). 
These results are the same as in the case
of the nonlinear Schr\"odinger equation on the real line. In fact, as regards the existence of  ground states, the results for systems on graphs differ, in general, from the ones for systems on the line even in the subcritical case: in the latter case, whenever the
constrained energy is lower bounded there always exist ground states (the {\em solitons}, whose shape is explicitly known), 
whereas for graphs the existence of a ground state is not guaranteed.

More precisely, we show that the existence of such constrained ground states is strongly conditioned by the topology of the graph. In particular, 
in the subcritical case we single out a 
topological hypothesis that prevents a graph from having ground states for every value of the mass. 

For the critical case, our results show a phenomenology much richer than the analogous on the line:  if some topological assumptions are fulfilled,
then there may exist a whole interval of masses for which a ground state exist. This behaviour is highly non-standard for $L^2$-critical 
nonlinearities.

\medskip

\medskip

\noindent \textbf{Mathematics~Subject~Classification~(2010):}
35R02, 35Q55, 81Q35, 49J40.

}
\end{minipage}
\end{center}

\bigskip






\section{Introduction: why dynamics on networks?}
Evolution on {\em metric graphs}, or {\em networks},
is a mathematical model used in order to approximate the dynamics
of systems located on {\em branched} spatial structures.
Such structures are characterized by the fact that locally
only one direction is important, except for some points
 where 
several directions  are available. Such special points in the structure are called {\em vertices, nodes} or {\em bifurcation points}, and the (possible) connections between two of them are called {\em edges}.

The research on dynamics of networks
started in 1953 with the seminal work by Ruedenberg and Scherr (\cite{RS1953}) where  the
dynamics of valence electrons in organic molecules was approximated by defining a  suitable Schr\"odinger
operator on the molecular bonds, treated as edges of a metric graph. This paper initiated the
research line nowadays known as {\em evolution on quantum graphs} (see the milestone paper by
Kostrykin and Schrader \cite{KS99} and the treatise by Berkolaiko and Kuchment \cite{berkolaiko}).
By definition, a quantum graph is a network, made of edges and vertices, on which functions are defined
and a linear differential operator acts.

More recently, several papers appeared, in which a nonlinear evolution on branched structure was proposed 
(\cite{bona,acfn11,acfn12,acfn14,acfn16,cfn,ps,nps,p,marzuola,matrasulov,sobirov,caudrelier,tentarelli,tentarelliserra,sven}). The first systematic
study of nonlinear dynamics on networks is contained in \cite{alimehmeti}, however it is only in the last
two years that a great deal of efforts in this direction has been carried out.

Here we focus on the problem of establishing the existence of ground states for the Nonlinear Schr\"odinger (NLS) Equation on metric graphs.
In particular, we shall review the results given in \cite{ast15-2,ast16-1,ast16-2}. 

The problem generalizes along two directions
the issue of finding the ground state of a Bose-Einstein condensate: first, the chosen domain is not standard, as it consists of a network instead of a three-dimensional regular region, or a disc, or a "cigar''; second, the nonlinearity we consider in the energy functional displays an arbitrary power, whilst
typically for condensates in the so-called Gross-Pitaevskii regime the quartic power emerges as the effective one. Furthermore, we limit our analysis to the
 {\em focusing} case, i.e. the case in which the net effect of the nonlinearity on the
time evolution is the concentration of the wave packets. 

The present note is organized as follows: in the rest of the introduction we 
give a historical overview on the mean-field limit for a many-boson system and 
draw a link between the problem of minimizing the constrained energy, possibly on graphs, and the Bose-Einstein 
condensation; we finally give the basic definitions and notation and state the problem we shall focus on. In Section 2 we
illustrate in a rather formal way the role of the so-called {\em critical nonlinearity power}, then we write down the Euler-Lagrange equation and
the Kirchhoff condition. In Section 3 we give some well-known and some others less known examples in which the problem is solved, trying to convey 
some general ideas. Section 4 is devoted to the introduction of a topological hypothesis (Assumption (H)) that is the core of the key nonexistence
result given in Theorem \ref{nonex}, to which
Section 5 is devoted together with a review of rearrangement theory (for non-experts). In Section 6 we give many examples in which Assumption (H) is not satisfied
and show how to prove existence of ground states through a technique of {\em graph surgery}. The point we stress here is that in all cases where topological information is not sufficient to solve the problem,
analysis needs to be carried out case by case; to this aim, Theorem 6.2 with its operative Corollary 6.1 can give a great help, as it states that in order to ensure the existence of a ground state it is sufficient to find a state that {\em does better than a soliton on the line}, i.e. a state whose  energy level is  lower than the level of the solution to the same problem on the line.  Finally, in Section 7 we treat the case of the critical power 
nonlinearity, where the role of topology overwhelms that of the metric, at least for simple cases, but the lack of compactness of minimizing sequences is 
much more serious. The analysis becomes thus more involved but, as a result, graphs can be classified in four disjoint categories, for each of those we give an
exhaustive result (Theorems 7.1-7.4).

\subsection{Nonlinearity and Condensation}

It is nowadays well-established both theoretically (\cite{bose,einstein,stringari}) and
experimentally (\cite{cornell,ketterle}) that, at a critical (usually very low, amounting to few Kelvin)
temperature, an ultracold gas of identical bosons (e.g. atoms and ions like sodium, rubidium and potassium) in a magnetic and/or optical trap experiences a phase
transition that turns the system into a {\em Bose-Einstein condensate}, i.e. a phase in which
a macroscopic fraction of the elementary components acquires a one-particle quantum state (i.e. a {\em wave function
${\varphi}$}): furthermore, the bosonic symmetry imposes that such a state {\em is the same
for all particles}. The system then can be thought of as a unique {\em giant quantum particle} lying
in the state $\varphi$, called {\em ground state of the condensate}. Such a state can be found as a solution to the
{variational problem}
$$
\min_{u \in H^1 (\Omega), \int_\Omega u^2 = N} E_{GP} (u)
$$
where the {\em Gross-Pitaevskii functional} $E_{GP}$ reads
\begin{equation}
    \label{gpfunctional}
E_{GP} (u) \ = \ \f 1 2 \int_\Omega | \nabla u (x) |^2 \, dx + 8 \pi \alpha \int_{\Omega} | u (x) |^4 \, dx.
\end{equation}
Here $ \alpha$ is the {\em scattering length} of the two-body interaction between the particles of the gas,
$\Omega$ is the spatial domain defined by the trap, and $N$ is the number of particles in the condensate.

Sometimes, the presence of the
trap is not modeled just by bounding the integral to the domain $\Omega$, but rather endowing the functional with an additional term that takes into account the presence of a confining potential, often
a harmonic one.

The rigorous derivation of the
 functional \eqref{gpfunctional} is the core result of the {\em Gross-Pitaevskii theory}, that is an {\em effective} theory used in order to describe the behaviour of a Bose-Einstein condensate. As we shall point out later in some more detail, the main merit of such a theory is that it 
reduces the complexity of the problem from $N$-body to one-body, even though the 
resulting system is nonlinear.

In the first experimental realizations of condensation, the shape of the trap was definitely three-dimensional
and regular. Since then, the technology of traps underwent
an impressive development, so that nowadays disc-shaped and cigar-shaped traps are currently produced, and
some indication of the occurrence of a Bose-Einstein condensation on a ramified structure (i.e. in
a Josephson junction) has been recently provided
(\cite{lorenzo}).

\subsection{From a linear $N$-body to a nonlinear 1-body problem}

The dynamics of the Bose-Einstein condensates naturally raises a question: given that quantum mechanics is a {\em linear} theory, where does the nonlinearity (i.e. the quartic power in the energy functional \eqref{gpfunctional}) come from?

 Such a  problem can be formulated in the time-independent or in the time-dependent framework.

The time-independent formulation is closer to the problem  of the ground state. The validity
of the Gross-Pitaevskii theory and of the functional \eqref{gpfunctional} was rigorously established in a series of works by E.-H. Lieb,
R. Seiringer and J. Yngvason (see e.g. \cite{ls,lsy,ly}). Among their 
achievements, we recall that they found the correct scaling for the potential
describing the pair interaction between the particles, namely
\begin{equation} \label{scaling}
V (x_i - x_j) \ \longrightarrow V_N (x_i - x_j) \ : = \ N^2 V (N (x_i - x_j)) ,
\end{equation}
so that the scattering length of the interaction scales as $1/N$.

A celebrated result in the cited paper is the proof of the Bose-Einstein condensation {\em
in the ground state} of the $N$-body system of bosons. This means that in the
ground state of the $N$-body
Hamiltonian operator
\begin{equation} \nonumber 
 H_N \ = \ \sum_{j=1}^N ( - \Delta_{x_j} + W (x_j) ) + \sum_{i < j} V_N (x_i - x_j)
\end{equation}
representing the energy of the boson gas, 
the $k$-particle correlation function converges to the factorized state $\varphi(x_1) \dots
\varphi (x_k)$ as the number of particles $N$ grows to infinity.
The resulting function $\varphi$ minimizes the constrained functional \eqref{gpfunctional}.

In physical terms,
in the limit $N \to \infty$ all particles collapse in the same quantum state,
represented by a wave function $\varphi$ minimizing the Gross-Pitaevskii energy functional \eqref{gpfunctional}.
Thus, one has {\em condensation in the ground state}.

On the other hand, the emergence of factorized states and of a nonlinearity out of a linear dynamics can be described also in the time-dependent
framework. The time-dependent formulation of the problem of the description of a dilute boson gas
in the Gross-Pitaevskii regime historically arises
from
one of the most topical and active fields of research of the mathematical physics in the last two decades, namely the
{\em mean field limit} for the dynamics  of many-body systems. The problem can be
summarized as follows: one is interested
in studying the evolution in time of a system made of a huge number $N$ of identical particles. According to the
basics of quantum mechanics, the state of the whole system 
at time $t$ is represented by a wave function
$\Psi_N (t,x_1,\dots, x_N)$, where $x_i$ is the position variable of the $i.$th particle and the
evolution of $\Psi_N$ is described by the $N$-body, linear Schr\"odinger equation
\begin{equation} \begin{split}
\label{N-schrod}
i \partial_t \Psi_N(t,x_1,\dots, x_N) \ = \  & -\sum_{j=1}^N \Delta_{x_j}  \Psi_N(t,x_1,\dots, x_N)
\\ &
+ \sum_{i < j} V (x_i - x_j)\Psi_N(t,x_1,\dots, x_N) \end{split}
\end{equation}
where the potential $V$ models the (two-body) interaction between the particles.

Such an equation is in general impossible to solve or even to study numerically, due to the fact that
$N$ is often very large (from around millions to $10^{23}$).
However, it is well-known that physicists are used to deal with such a systems by reducing the equation from $N$-body to one-body, but paying the price of introducing a nonlinear term in the equation. Justifying
such an approximation and providing an estimate for the error made when employing it, has been the main task of the research on mean-field limit for large systems of identical interacting bosons. In its
simplest version, the problem of the mean-field limit can be expressed as follows: prove that the evolution provided
by the equation \eqref{N-schrod} with the potential modified through a {\em weak coupling scaling} $V \longrightarrow V/N$, i.e.
\begin{equation}  \nonumber
\begin{split}
i \partial_t \Psi_N(t,x_1,\dots, x_N) \ =  & \ -\sum_{j=1}^N \Delta_{x_j}  \Psi_N(t,x_1,\dots, x_N)
\\ & \
+ \frac 1 N \sum_{i < j} V (x_i - x_j)\Psi_N(t,x_1,\dots, x_N)
\end{split}
\end{equation}
with the {\em factorized} initial data
\begin{equation} \nonumber
\Psi_0 (x_1, \dots. x_N) \ = \ \phi_0 (x_1) \dots \phi_0 (x_N)
\end{equation}
can be approximated by the evolution of $N$ independent particles, following each the dynamics given by
\begin{equation} \label{hartree}
i \partial_t \phi (t,x) = - \Delta \phi (t,x) + (V \star |\phi (t,x)|^2) \phi (t,x).
\end{equation}
More precisely, the task is to prove that in the limit $N \to \infty$ the correlation functions of the $k$-particle subsystems converge to
$   {\phi (t,x_1) \dots \phi (t,x_k)}$ where $   {\phi (t,x)}$ solves \eqref{hartree}.
The history of the main achievements of this research line starts with the work of
 K. Hepp (\cite{hepp}), who stated the problem (even though his celebrated work is rather devoted to the classical limit of quantum mechanics). Then,
  J. Ginibre and G. Velo (\cite{ginibrevelo79}) treated the problem of mean field for Coulombian systems in the second-quantized framework. On the other hand,
 H. Spohn (\cite{spohn80}) proved the mean-field limit for systems of particles interacting through a bounded potential, using a first quantization formalism.
   In 2000, C. Bardos, F. Golse and N. Mauser (\cite{bgm00}) re-obtained Spohn's result by splitting the problem in the issue of the convergence as $N$ goes to 
infinity of the $N$-body Schr\"odinger equation to an infinite hierarchy, and in the problem of the uniqueness of the solution to such resulting hierarchy. For Coulombian systems, 
the latter problem was solved one year later by
    L. Erd\H{o}s and H.-T. Yau (\cite{ey00}). 

All the cited works are purely mean-field, so that the main result they get is equation \eqref{hartree}, that bears a non-local nonlinearity.
The first result that opened the road to the derivation of an effective equation with a local nonlinearity was given
in \cite{eesy}, where a mixed scaling  $V \longrightarrow N^{3\gamma} V(N^\gamma \cdot) / N$ was introduced. 
This can be {\em formally} interpreted as a 
mean-field theory for a 
smoothened Dirac's delta. The effective equation yielded by the limit is then \eqref{hartree} with the replacement of $V$ with the Dirac's delta
potential, so that
    \begin{equation} \nonumber
i \partial_t \phi (t,x) = - \Delta \phi (t,x)
+ \left( \int V dx \right)|\phi (t,x)|^2 \phi (t,x)
\end{equation}
is the target dynamics: this is the time-dependent
Gross-Pitaevskii equation that describes the evolution of Bose-Einstein condensates.
 The derivation was not complete, since the problem of the uniqueness of the solution to the 
limit hierarchy, stated in (\cite{bgm00}) was not solved.

\noindent
Later,
 R.A., F. Golse and A. Teta (\cite{agt})
derived the cubic Schr\"odinger equation in dimension one by
studying a scaling limit for a system of one-dimensional identical bosons interacting through a repulsive pointwise interaction: morally, again, the limit whose existence they proved is a mean-field limit with
    a Dirac's delta potential. Since the result is limited to one-dimensional systems, its validity is restricted to  {\em cigar-shaped} condensates. A more general result, valid
    for attractive interaction too, was given later in \cite{holmer}.

\noindent
Finally, in a series of works dating from 2006 to 2010
(\cite{esy1,esy2,esy3}), L. Erd\H{o}s, B. Schlein and H.-T. Yau derived the Gross-Pitaevskii equation for three-dimensional systems. The scaling they adopted is the one discovered by E.-H. Lieb and R. Seiringer for the time-independent framework \eqref{scaling}.

Once derived the effective equation, it remained to estimate the error made by replacing the $N$-body linear Schr\"odinger equation by the one-body nonlinear
equation: for the pure mean-field scaling, the breakthrough came by
 I. Rodnianski and B. Schlein \cite{rs09}, who provided a new proof of the mean field limit, inspired to the work of Ginibre and Velo \cite{ginibrevelo79} and gave also the first estimate of the error.
In 2010, A. Knowles and P. Pickl \cite{kp10} gave a further proof of the limit, in the first-quantized formalism but avoiding use of hierarchies. The method allowed to deal with more singular potentials and produced a new estimate of the error.

\noindent
Further improvements on the rate of convergence for the mean field have been achieved in \cite{falconi, ammarinier,gmm1,gmm2,bucholz,nam1}, while for the Gross-Pitaevskii regime 
a similar estimate has been proved in \cite{bos}. For a complete review on most recent results see the monography \cite{schleinlibro}

\medskip

In this review there are no proofs of the condensation or of the Gross-Pitaevskii regime for systems on graphs.
In fact,
it is widely known that no condensation can occur in one-dimensional systems, in the sense that the phase transition that defines the condensation cannot take place: however, in 2015
J. Bolte and Kerner \cite{boltekerner} proved condensation for {\em free} gas and no condensation
for interacting gases in graphs, considered as quasi-dimensional systems, in the sense of the presence of the phase transition. Moreover, in 1996 W. Ketterle and N.J. van Druten \cite{kvd} gave an evidence of a concentration phenomenon under some aspects analogous to condensation. Finally, one-dimensional condensates can be considered as squeezing limits of three-dimensional condensates, as proved by R. Seiringer and Lin \cite{seiringerlin}.

\subsection{The problem}

In these lectures we search for the simplest solutions to the {\em nonlinear focusing} Schr\"odinger equation on graphs, avoiding the problem of the rigorous derivation from first principles. Simplest solutions are particular cases of {\em standing waves}, that minimize the energy under some physical constraints:
in particular, we consider the so-called {\em mass constraint}.

Before stating the problem precisely, we need some definitions and notation.

\begin{figure}[t]
\begin{center}
\begin{tikzpicture}
\node at (0,0) [nodo] (1) {};
\node at (-1.5,0) [infinito]  (2){$\infty$};
\node at (1,0) [nodo] (3) {};
\node at (0,1) [nodo] (4) {};
\node at (-1.5,1) [infinito] (5) {$\infty$};
\node at (2,0) [nodo] (6) {};
\node at (3,0) [nodo] (7) {};
\node at (2,1) [nodo] (8) {};
\node at (3,1) [nodo] (9) {};
\node at (4.5,0) [infinito] (10) {$\infty$};
\node at (5.5,0) [infinito] (11) {$\infty$};
\node at (4.5,1) [infinito] (12) {$\infty$};
\node at (1.4,0.7) [nodino]  {};

\draw [-] (1) -- (2) ;
\draw [-] (1) -- (3);
\draw [-] (1) -- (4);
\draw [-] (3) -- (4);
\draw [-] (5) -- (4);
\draw [-] (3) -- (6);
\draw [-] (6) -- (7);
\draw [-] (6) to [out=-40,in=-140] (7);
\draw [-] (3) to [out=10,in=-35] (1.4,0.7);
\draw [-] (1.4,0.7) to [out=145,in=100] (3);
\draw [-] (6) to [out=40,in=140] (7);
\draw [-] (6) -- (8);
\draw [-] (6) to [out=130,in=-130] (8);
\draw [-] (7) -- (8);
\draw [-] (8) -- (9);
\draw [-] (7) -- (9);
\draw [-] (9) -- (12);
\draw [-] (7) -- (10);
\draw [-] (7) to [out=40,in=140] (11);
\end{tikzpicture}
\end{center}
\caption{How a connected, non-compact metric graph may look like.}
\label{how}
\end{figure}
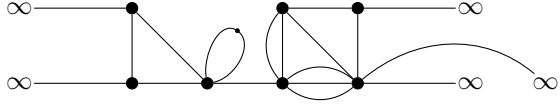

\begin{itemize}
    \item A {\em{graph}} $\G$ i.e. a couple of sets $(\mathcal V,\mathcal B)$, where
    $\mathcal B$ is a subset of $\mathcal V \times \mathcal V$.

    The set $\mathcal V$ is interpreted as the set of the {\em vertices}, i.e. points in the space, while $\mathcal B$ is the set of the {\em edges} or {\em bonds}, e.g. links between vertices: every edge is
    then identified with the couple of vertices it connects.

    \noindent The {\em degree} of  a vertex $\textsc v \in \mathcal V$ is the number of edges
    starting from or ending at $\textsc v$.

    \noindent
    Both $\mathcal V$ and $\mathcal B$ are {\em finite} sets, so that we shall always deal with
    graphs with a finite number of edges and vertices.

    \noindent
    There is a non-empty subset $\mathcal V_\infty$ of $\mathcal V$, made of {\em vertices at infinity}. Two vertices at infinity cannot be connected by edges, and every vertex at infinity has degree equal to one.
    We shall refer to edges ending at a vertex at infinity as to {\em halflines}.

    \noindent
    Two graphs are topologically equivalent if they can be deformed into
each other without changing the sets {$\mathcal V, \mathcal V_\infty$} and {$\mathcal B$}.

    \item
    In order to construct a {\em metric graph}, an edge $e$ is identified with an interval $I_e : = [0, \ell_e]$, where $\ell_e \in [0, + \infty]$. This
    correspondence fixes the {\em metric} of the graph.

    \noindent{Given a topology, several metrics are possible, as every
finite edge can have an arbitrary length. Conversely, the metric on the
halflines is fixed.

\noindent{For a pictorial idea of a generic graph see Fig. \ref{how}, that show an example where selfloops, multiple connections, and haflines are present.}
}

    \item
    A function ${u: \G \to  {\mathbb{C}}}$ is a bunch of functions ${u  = (u_e)_{e \in \mathcal B}}$, with
$
u_e : I_e \to  {\mathbb{C}} 
$.

\noindent
A function $u$ is continuous on $\G$ if every $u_e$ is continuous in $I_e$ and if $u$ is continuous at vertices, namely, if the value attained at a vertex $\textsc v$ is independent of the edge chosen to reach $\textsc v$.

\item
A function $u : \G \to \mathbb C$ is integrable if every function $u_e$ is integrable on $I_e$, and
$$ \int_\G u \, dx \ : = \ \sum_{e \in \mathcal B} \int_{0}^{\ell_e} u_e (x_e) \, dx_e
$$
The usual functional spaces can be defined as
$$
L^p (\mathcal G) \ := \ \bigoplus_{e \in \mathcal B} L^p (I_e), \qquad
 \|u\|_{L^p(\GG)}^p  := \sum_{e\in \mathcal B}\|u_e\|_{L^p(I_e)}^p.
$$
The space
${H^1 (\mathcal G)}$ is defined as the set of {continuous} functions ${u  = (u_e)_{e \in \mathcal B}}$  such that
$$
{u_e \in H^1(I_e)} \quad\forall e \in \mathcal B,\qquad\quad
\|u\|_{H^1(\GG)}^2= \sum_{e\in \mathcal B}\|u_e\|_{H^1(I_e)}^2.
$$
We stress that
continuity is imposed at vertices too, so that no jump can occur.

\item Fixed $\mu > 0$, we define
    $$
    H^1_\mu (\G) : = \{ u \in H^1 (\G), \, \|u \|^2_{L^2 (\G)} = \mu \},
    $$
    that is the space of the functions in $H^1 (\G)$ that fulfil the {\em mass constraint}.
   \end{itemize}

\noindent
On the {metric graph $\mathcal G$} let us define
$$ E (u, {\mathcal G}) = \f 1 2 \| u' \|_{L^2 (\mathcal G)}^2 - \f 1 p
\| u \|_{L^p (\mathcal G)}^p, 
$$
that is a functional in  ${C^1(H^1_\mu(\G),\R)}$ for all ${p\in [2,+\infty)}$.

\medskip

\noindent
The problem we treat is the following:

\medskip

\noindent{\bf{Problem P.}}
Given a {\em connected, non-compact} metric graph $\G$ and fixed $\mu > 0$,
does there exist a
{\em ground state at mass $\mu$}, namely
a
minimizer of {$E (\cdot, \G)$} in the space $H^1_\mu (\G)$?

\medskip

\noindent
In other words, we look for functions ${ u\in  H^1_\mu(\G)}$ such that

$$
E(u,\G) \ =  \ \E_\G (\mu) 
$$

\noindent{where we introduced the notation} 
\begin{equation} \label{superE}
{\E_\G (\mu) \ : = \ \inf_{v\in H^1_\mu(\G)} E(v,\G)}.
\end{equation}

Notice that, as regards the problem of the existence of a minimizer,
if $\G$ is compact, then every minimizing sequence is compact, so that a minimizer always exists. For this reason we restrict to non-compact graphs, i.e. graphs that contain at least one halfline. On the other hand, it will be clear from the analysis that if a graph is not connected, then minimizers concentrate on the most convenient connected component, so that the hypothesis that the graph is connected is not restrictive.

We end this introductory part by noticing that, due to the shape of the energy functional, we shall limit
the analysis to nonnegative functions.

%




%

\section{Preliminary remarks}

\subsection{Mass constraint and lower boundedness}

First of all observe that, regardless of the chosen graph
$\G$,  if no constraint is imposed then the functional $E (\cdot, \G)$ is {\em not} lower bounded. Indeed, fixed $u \in H^1 (\G)$,
one has
$$ E (\lambda u, \G) = \f {\lambda^2} 2 \| u' \|^2_{L^2 (\G)} - \f { \lambda^p} p
\| u \|^p_{L^p (\G)} \ \to \ - \infty, \quad \lambda \to + \infty.$$

\noindent
On the other hand,
let us  restrict to the case of star graphs made of halflines, and impose the mass constraint $\| u \|_{L^2 (\G)}^2 = \mu$. On these particular graphs,
one can perform mass preserving
transformations ${u(x) \to \sqrt \lambda u (\lambda x) : = u_\lambda (x)}$, so that
$$ E (u_\lambda, \G) = \f {\lambda^2} 2 \| u' \|^2_{L^2 (\G)} - \f { \lambda^{\f p 2 -1}} p
\| u \|^p_{L^p (\G)}. $$

\noindent
Now,
\begin{itemize}
\item
if $  {p < 6}$, then kinetic energy prevails and this suggests that the  energy is {\em lower bounded}.
Indeed by using Gagliardo-Nirenberg estimates one can prove that this is actually the case. The problem
with $2 < p < 6$ is referred to as the {\em subcritical case}.

\item
If $  {p > 6}$, then the potential term overwhelms the kinetic one and $ E (u_\lambda, \G) \to - \infty$,
as $\lambda \to + \infty$, so the energy is not lower bounded. The problem
with $p > 6$ is then referred to as the {\em supercritical case}.

\item
If $  {p = 6}$, then there is a delicate balance between kinetic and nonlinear term. As we shall see, lower boundedness
of $E$ depends on the value of $\mu$. The problem
with $p = 6$ is referred to as the {\em critical case}.
\end{itemize}

\subsection{The Euler-Lagrange equation: Kirchhoff's rule}
As a minimum of the {\em constrained} functional,
every ground state $   u$ must satisfy the Lagrange Multiplier
Theorem
$$ \nabla E (u, \G) (u) = \frac \omega 2 \nabla (\mu - \| u \|^2_{L^2 (\G)} ) $$
for some ${\omega \in \R}$ (notice that in order to simplify notation here we called $\frac \omega 2$ 
the Lagrange multiplier). Then, for every ${\eta \in H^1 (\G)}$,
\begin{equation*}
\begin{split}
0 \ = \ &    \nabla E (u, \G) \eta - \frac \omega 2 \nabla (\mu - \| u \|^2) \eta
\  = \ \int_\G ( u' \eta' - u^{p-1} \eta + \omega u \eta) \, dx \\
 \ = \ & \sum_{e \in \mathcal B}  \int_0^{\ell_e} ( u_e' \eta_e' - u_e^{p-1}  \eta_e + \omega u_e
\eta_e) \, dx_e \\
\ = \  &\sum_{e \in \mathcal B} \ u'_e \eta_e |_0^{\ell_e} 
+ \sum_{e \in \mathcal B} \int_0^{\ell_e} ( - u_e'' - u_e^{p-1} + \omega u_e ) \eta_e \, dx_e
        \end{split}
\end{equation*}
Notice that the first term concerns vertices, while the second is determined by the values of the
integrand inside the edges.

Now pick an edge $  {\bar e}$ and consider a function $  {\eta \in C_0^\infty (\bar e)}$.
Then the second term only survives and forces the
{\em Stationary Nonlinear Schr\"odinger equation}
\begin{equation} \label{el}
u_{\bar e}'' + u_{\bar e}^{p-1} = \omega u_{\bar e}
\end{equation}
{to hold inside $\bar e$, and, by arbitrarity of $\bar e$, on every edge}.

Now consider $   {\eta \in C^\infty_0 (\G)}$ that vanishes {\em in all vertices except one},
say $  {\bar \vv}$. Then, boundary terms survive {\em if and only if they refer to edges starting from or
ending at $\bar \vv$.} Their contribution reads
\begin{equation*}
    \begin{split}
  &  \sum_{e \in E} \ u'_e \eta_e |_0^{\ell_e} =  \sum_{e \rightarrow \bar \vv}
u_e' (\ell_e) \eta (\ell_e) -   \sum_{e \leftarrow \bar \vv}
u_e' (0) \eta (0) \\ & \
=  \ \left(  \sum_{e \rightarrow \bar \vv}
u_e' (\ell_e)  -   \sum_{e \leftarrow \bar \vv}
u_e' (0)
\right) \eta (\bar \vv)
    \end{split}
\end{equation*}
where we denoted $e \rightarrow \bar \vv$ the edges ending at $\bar \vv$ and $ e \leftarrow \bar \vv$
the edges starting from $\bar \vv$.

\noindent
Thus
$$\sum_{e \rightarrow \bar \vv}
u_e' (\ell_e)  -   \sum_{e \leftarrow \bar \vv} u_e' (0) \ = \ 0$$
that is called Kirchhoff's rule, and is often expressed by saying that at every vertex the global ingoing (or outgoing) derivative
vanishes. A common compact form of Kirchhoff's rule is
\begin{equation} \label{kirch}
\sum_{e\succ \vv}
\frac{d u_e}{d x_e}(\vv)=0.
\end{equation}

\medskip \noindent
Let us finally recall that the Euler-Lagrange equations \eqref{el} together with the Kirchhoff's conditions
can be summarized in the equation
\begin{equation} \label{lapkir}
\Delta_K u (x) + u^{p-1} (x) = \omega u (x)
\end{equation}
where $\Delta_K$ is the operator acting as the laplacian on functions that are $H^2$ on every edge and fulfil Kirchhoff's rule at
all vertices. It is immediately seen that equation \eqref{lapkir} is the stationary equation associated to the time-dependent
nonlinear Schr\"odinger equation
\begin{equation} \label{schrod}
i \partial_t \psi (x) \ = \ - \Delta_K \psi (t) - |\psi (t)|^{p-2} \psi (t),
\end{equation}
and it is well-known that for such equation the dynamics preserves the $L^2$-norm and the value of the energy $E(\cdot, \G)$.

We stress that \eqref{lapkir} is equivalent to the condition of stationarity of the functional, so that it is satisfied not only by  ground states,
but also by every standing wave of the nonlinear Schr\"odinger equation, namely by all solutions
to \eqref{schrod} of the type
\begin{equation} \nonumber 
 \psi (t,x) = e^{i \omega t} u(x),
\end{equation}
so that it is clear that the Lagrange multiplier $\omega$ has the dynamical meaning of a {\em frequency}.

\section{Examples}
In this section we give some basic examples of graphs and the related results concerning the
existence or the nonexistence of ground states.

\subsection{The real line}

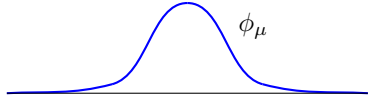
\begin{figure} 
\begin{center}
\begin{tikzpicture}[xscale= 0.4,yscale=0.6]
\draw[-,thick,blue] (0,2) to [out=0,in=170] (2.5,0.2);
\draw[-,thick,blue] (0,2) to [out=180,in=10] (-2.5,0.2);
\draw[-,thick,blue] (2.5,0.2) to [out=350,in=178] (6,0);
\draw[-,thick,blue] (-2.5,0.2) to [out=190,in=2] (-6,0);
\draw[-] (-6,0) -- (6,0);
\node at (2.2,1.5) {$\phi_\mu$};
\end{tikzpicture}
\end{center}
\caption{The soliton $\phi_\mu$.} \label{fig:soliton}
\end{figure}

It is well-known (\cite{zs72,cl82,gss1}) that
for ${p\in (2,6)}$ and ${\mu>0}$ {ground states exist} and are
all translated of the {soliton} (Fig. \ref{fig:soliton})
$$
\phi_\mu (x) \ = \ C \mu^{\f 2 {6-p}} \hbox{\mbox{sech}}^{\f 2 {p-2}} (c
\mu^{\f {p-2} {6-p}} x).
$$
where $C$ and $c$ are irrelevant constants dependent on $p$ only and not on $\mu$.

\begin{figure}
\begin{tikzpicture}[xscale=0.7,yscale = 0.7]
\begin{axis}[axis lines=middle,enlargelimits,xtick={-6,6},
ytick={-40,10}]
\addplot
[domain=-5.2:5.2,samples=40,smooth,ultra thick,blue]
{x^4 /4 - 13 * x^2 / 2 };
\addplot
[domain=-5.2:5.2,variable=\t,
 samples=40,smooth,ultra thick,brown]
({t},5);
\addplot
[domain=-4.8:-1.8,variable=\t,
 samples=40,smooth,ultra thick,darkgray]
({t},-20);
\addplot
[domain=4.8:1.8,variable=\t,
 samples=40,smooth,ultra thick,green]
({t},-20);
\addplot
[domain=-5:5,variable=\t,
 samples=40,smooth,ultra thick,red]
({t},0);
\end{axis}
\end{tikzpicture}
\qquad
\begin{tikzpicture}[xscale=0.7,yscale = 0.7]
\begin{axis}[axis lines=middle,enlargelimits,xtick={-6,6}]
\addplot
[domain=0:6,variable=\t,
 samples=500,smooth,ultra thick,red]
({t},{sqrt(13 * t^2 / 2- t^4 /4)});
\addplot
[domain=0:5.3,variable=\t,
 samples=500,smooth,ultra thick,red]
({-t},{sqrt(13 * t^2 / 2- t^4 /4)});
\addplot
[domain=0:5.3,variable=\t,
 samples=500,smooth,ultra thick,red]
({t},{-sqrt( 13 * t^2 / 2 - t^4 /4)});
\addplot
[domain=0:5.3,variable=\t,
 samples=500,smooth,ultra thick,red]
({-t},{-sqrt( 13 * t^2 / 2 - t^4 /4)});
\addplot
[domain=-5.5:5.5,variable=\t,
 samples=1000,smooth,ultra thick,brown]
({t},{sqrt(12.5 + 13 * t^2 / 2- t^4 /4)});
\addplot
[domain=-5.5:5.5,variable=\t,
 samples=1000,smooth,ultra thick,brown]
({t},{-sqrt(12.5 + 13 * t^2 / 2- t^4 /4)});
\addplot
[domain=-4.8:-1.8,variable=\t,
 samples=400,smooth,ultra thick,darkgray]
({t},{-sqrt(-18.44 + 13 * t^2 / 2- t^4 /4)});
\addplot
[domain=-4.8:-1.8,variable=\t,
 samples=400,smooth,ultra thick,darkgray]
({t},{sqrt(-18.44 + 13 * t^2 / 2- t^4 /4)});
\addplot
[domain=1.8:4.8,variable=\t,
 samples=400,smooth,ultra thick,green]
({t},{-sqrt(-18.44 + 13 * t^2 / 2- t^4 /4)});
\addplot
[domain=1.8:4.8,variable=\t,
 samples=400,smooth,ultra thick,green]
({t},{sqrt(-18.44 + 13 * t^2 / 2- t^4 /4)});
\end{axis}
\end{tikzpicture}
\caption{Left: The profile of the potential $V$. Right: The phase portrait induced by the potential $V$.}
\label{portrait}
\end{figure}
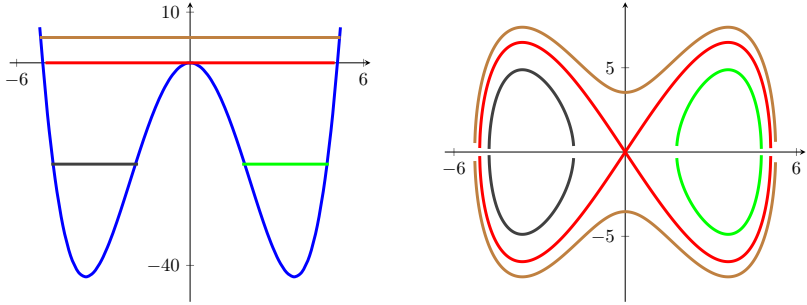

\noindent
If ${p=4}$,  i.e. in the case of the cubic Schr\"odinger equation, one gets
$$
\phi_\mu (x) \ = \frac{\mu}{2\sqrt2} \hbox{\mbox{sech}}\left(\frac{\mu}4
x\right), \qquad E (\phi_\mu, \R) = -\frac{\mu^3}{96}.
$$
Let us point out that the solitons and their translated are the only stationary solutions
to \eqref{schrod} on the line. In order to prove it, notice that
the Euler-Lagrange equation on the line
$${ u'' + u^{p-1} = \omega u},$$
can be rewritten as
\begin{equation} \label{newton}
u'' = - \frac {dV}{du} (u)
\end{equation}
where we defined
$$ V (u) = \f {u^p} p - \omega \f {u^2} 2 $$
Equation \eqref{newton} can be interpreted as a mechanical conservative problem, whose
phase portrait is displayed in Fig. \ref{portrait}.


\noindent
It is then clear that all solutions are periodic (and therefore not in $H^1(\R)$)
except the non-constant ones contained in the
{\em separatrix}, corresponding to solutions to \eqref{newton} with vanishing mechanical
energy. They turn out to be
the solitons and their translated.




\subsection{The halfline}

\begin{figure} 
\begin{center}
\begin{tikzpicture}[xscale = 0.5, yscale=0.8]
\draw[-,thick,blue] (0,2) to [out=0,in=170] (2.5,0.2);
\draw[-,thick,blue] (2.5,0.2) to [out=350,in=178] (6,0);
\draw[-] (0,0) -- (6,0);
\node at (2.4,1.3) {$\phi_{2\mu}$};
\end{tikzpicture}
\end{center}
\caption{Half a soliton.}
\label{fig:halfsoliton}
\end{figure}
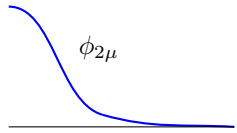

In the case $\G = \R^+$,
 ${p\in (2,6)}$ and ${\mu>0}$, there is exactly {\em one}  ground state
given by ``{half a soliton}'' (Fig. \ref{fig:halfsoliton}) of mass ${2\mu}$ (notice that in this case the translational symmetry is broken).

\noindent
If ${p=4}$, then
$$
\phi_{2\mu} (x) \ = \frac{\mu}{\sqrt2}
\hbox{\mbox{sech}}\left(\frac{\mu}2 x\right), \qquad E (\phi_{2\mu},
\R^+) = -\frac{\mu^3}{24}.
$$

\subsection{A star-graph made of halflines}

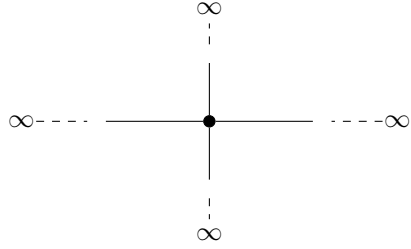
\begin{figure}[h]
\begin{center}
\begin{tikzpicture}
\node at (-2.5,0) [infinito]  (1) {$\infty$};
\node at (0,0) [nodo] (2) {};
\node at (2.5,0) [infinito]  (3) {$\infty$};
\node at (0,1.5) [infinito]  (4) {$\infty$};
\node at (0,-1.5) [infinito]  (8) {$\infty$};
\node at (1.5,0)  [minimum size=0pt] (5) {};
\node at (-1.5,0) [minimum size=0pt] (6) {};
\node at (0,0.9) [minimum size=0pt] (7) {};
\node at (0,-0.9) [minimum size=0pt] (9) {};

\draw[dashed] (1) -- (6);
\draw[dashed] (3) -- (5);
\draw[dashed] (7) -- (4);
\draw[dashed] (9) -- (8);
\draw [-] (2) -- (6) ;
\draw [-] (5) -- (2) ;
\draw [-] (2) -- (7) ;
\draw [-] (2) -- (9) ;
\end{tikzpicture}
\end{center}
\caption{A star-graph made of four halflines.}
\label{4star}
\end{figure}

The cases of the halfline and of the line naturally generalize to the case of the star-graph $\mathcal S_n $
made of $n$
halflines (the case $n=4$ is represented in Fig. \ref{4star}). Yet the result changes, as proved in \cite{acfn12,acfn14,acfn14-2,acfn16}.

\noindent
Indeed,
for ${p\in (2,6)}$ and ${\mu>0}$,
$$
\EE_{\mathcal S_n} (\mu) \ =  \ \EE_\R (\mu)
$$
but the infimum is {not} achieved, so that there is {no ground state}.

In order to prove it, restrict to the simplest case
$p=4$ and
consider the star-graph $\mathcal S_3$ made of three halflines $\mathcal H_1, \mathcal H_2$, and $\mathcal H_3$, and
a function $u \in H^1_\mu (\mathcal S_3)$.
Let us introduce the notation ${u = (u_1, u_2, u_3)}$, where $u_j$ is the restriction of $u$ to the halfline $\mathcal H_j$,
and
${\mu_j := \| u_j \|_{L^2 (H_j)}^2}$. With no loss of generality, suppose ${\mu_1 = \min (\mu_1, \mu_2, \mu_3)}$
and construct a function $\widetilde u \in H^1_\mu (\mathcal S_3)$ such that $E(\widetilde u, \mathcal S_3) \leq E (u, \mathcal S_3)$
proceeding as follows:

\medskip

\noindent
1. On $\mathcal H_1$, replace $u_1$ with the {half-soliton} $   {\chi_{\R^+} \phi_{2 \mu_1}}$

\medskip

\noindent
2. On $\mathcal H_2 \cup \mathcal H_3$, replace the couple of functions $(u_2,u_3)$ with the {soliton} ${\phi_{\mu_2+\mu_3}}$.

\noindent
At this point, on $\mathcal S_3$ set the function $  {(\chi_{\R^+} \phi_{\mu_1}, \chi_{\R^+} \phi_{\mu_2+\mu_3}, \chi_{\R^+}
\phi_{\mu_2+\mu_3}})$, that {may not be continuous}.

 \medskip

\noindent
3. {\em Translate} the soliton sat on $\mathcal H_1 \cup \mathcal H_2$ in order to obtain 
a continuous function $\widetilde u$ on $\mathcal S_3$.
It is immediate that $\widetilde u$ belongs to $H^1_\mu (\mathcal S_3)$.

\noindent
Then, exploiting the minimum properties of the half-soliton on the half-line, and of the soliton on the line,
\begin{equation*}\begin{split}
     E (u, \mathcal S_3) \ & = \ E (u_1, \mathcal H_1) + E ((u_2,u_3) , \mathcal H_2 \cup \mathcal H_3) \\
 & \geq \ E (\chi_{\R^+} \phi_{2 \mu_1}, \R^+) + E (\phi_{\mu_2 + \mu_3}, \R) \\
 & = \ - \frac {\mu_1^3} {24} - \frac {(\mu - \mu_1)^3} {96} \ = \ E (\widetilde u, \mathcal S_3).
\end{split}
\end{equation*}
By construction $  {\mu \geq 3 \mu_1}$, then the minimum is attained for
$$ \mu_1 \ = \ 0.$$
Therefore, minimizing sequences
concentrate on {\em one halfline only}, and
{\em reconstruct a soliton at infinity}. An analogous result can be obtained for
every $p < 6$.

This simple example provides at least two important messages: first, as regards the minimization of the energy,
it is not convenient to spread the wave functions on many edges.
Second, for every candidate ground state,
the soliton is a serious competitor!









%

\subsection{The 3-Bridge $\mathcal B_3$}

\begin{figure}[h]
\begin{center}
\begin{tikzpicture}
\node at (-3.5,0) [infinito]  (1) {$\infty$};
\node at (-0.2,0) [nodo] (2) {};
\node at (1.2,0) [nodo] (3) {};
\node at (4.5,0) [infinito]  (4) {$\infty$};
\node at (3,0)  [minimum size=0pt] (5) {};
\node at (-2,0) [minimum size=0pt] (6) {};
\draw[dashed] (6) -- (1);
\draw [-] (6) -- (2);
\draw [-] (5) -- (3) ;
\draw [-] (2) -- (3) ;
\draw[dashed] (5) -- (4);
\draw [-] (2) to [out=-40,in=-140] (3);
\draw [-] (2) to [out=40,in=140] (3);
\end{tikzpicture}
\end{center}
\caption{The 3-bridge $\mathcal B_3$.}
\label{3ponte}
\end{figure}
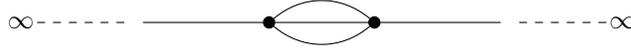

The first non-star cases treated in literature were the so-called {\em bridge graphs,}
first dealt with in \cite{ast15-1} (the triple bridge $\mathcal B_3$ is represented in Fig. \ref{3ponte}).
For $  {p\in (2,6)}$ and $  {\mu>0}$,
$$
\EE_{\mathcal B_3} (\mu) \ = \ \EE_\R (\mu)
$$
and again the infimum is   {not} achieved, so that there is   {no ground state}.

\medskip \noindent
\begin{figure}[h]
\begin{center}
\begin{tikzpicture}
\node at (-3.5,0) [infinito]  (1) {$\infty$};
\node at (-0.2,0) [nodo] (2) {};
\node at (1.2,0) [nodo] (3) {};
\node at (4.5,0) [infinito]  (4) {$\infty$};
\node at (3.2,0)  [minimum size=0pt] (5) {};
\node at (-2,0) [minimum size=0pt] (6) {};
\node at (-0.15,0) [minimum size=0pt] (7)  {};
\node at (1.15,-0.1)[minimum size=0pt] (8)  {};
\node at (1.15,0)[minimum size=0pt] (9)  {};
\node at (-0.15,0.1)[minimum size=0pt] (10)  {};

\draw (2) to [out=40,in=140] (8);
\draw [-] (7) -- (9);
\draw [-] (10) to [out=-40,in=-140] (3);
\draw [-] (6) -- (2);
\draw [-] (3) -- (5);

\draw[dashed] (6) -- (1);
\draw [dashed] (5) -- (4);
\end{tikzpicture}
\end{center}
\caption{Unfolding $\mathcal B_3$ into a line.}
\label{semieulerian}
\end{figure}
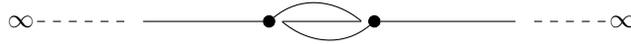

\noindent
Fixed $\mu > 0$
It turns out that for every function $u \in H^1_\mu (\mathcal B_3)$ one gets
$${E (u,  \mathcal B_3) \ > \ E (\phi_\mu, \erre)}. $$
To prove it, the key observation is that $\mathcal B_3$ is semi-Eulerian, namely, it can be unfolded into a line
together with every function $u \in H^1 (\mathcal B_3)$ (Fig. \ref{semieulerian}).

\noindent
So, let $u \in H^1_\mu (\mathcal B_3)$ and $\widetilde u \in H^1_\mu (\R)$ its unfolded version on
the line. Then, denoting by ${{\vv_1}, {\vv_2}}$ the two vertices,
on the line they correspond to points
$$x_1 \ < \  x_2 \ < \  y_1 \ < \ y_2$$
with $x_i, y_i$ associated to the vertex ${\vv_i}$.
Therefore, by continuity
$$ \widetilde u (x_1) = \widetilde u (y_1), \qquad
\widetilde u (x_2) = \widetilde u (y_2).
$$
Now,
if $  {\widetilde u (x_1) < \widetilde u (x_2)}$, then there is a minimum in the interval $(x_2, y_2)$.

\noindent
If $  {\widetilde u (x_1) > \widetilde u (x_2)}$, then there is a minimum in the interval $(x_1, y_1)$.

\noindent
If $  {\widetilde u (x_1) = \widetilde u (x_2)}$, then $u$ takes the same value in four points.

\smallskip \noindent
In all cases {$\widetilde u$ cannot be a soliton}, then
$$ E (u, \mathcal B_3) \ = \  E (\widetilde u, \R) \ > \ E (\phi_\mu, \R)  \ = \ \EE_\R  (\mu).$$
Nevertheless, it is possible to define a sequence that asymptotically reconstructs a soliton
on a halfline, so that
$$ \E_{\mathcal B_3} (\mu) \ = \ \E_\R (\mu) $$
but such a minimum is not attained.

\subsection{The double bridge $\mathcal B_2$}

\begin{figure}[h]
\begin{center}
\begin{tikzpicture}
\node at (-3.5,0) [infinito]  (1) {$\infty$};
\node at (-0.2,0) [nodo] (2) {};
\node at (1.2,0) [nodo] (3) {};
\node at (4.5,0) [infinito]  (4) {$\infty$};
\node at (3,0)  [minimum size=0pt] (5) {};
\node at (-2,0) [minimum size=0pt] (6) {};
\draw[dashed] (6) -- (1);
\draw [-] (6) -- (2);
\draw [-] (5) -- (3) ;
\draw[dashed] (5) -- (4);
\draw [-] (2) to [out=-40,in=-140] (3);
\draw [-] (2) to [out=40,in=140] (3);
\end{tikzpicture}
\end{center}
\caption{The two-bridge $\mathcal B_2$.}
\label{ponte2}
\end{figure}
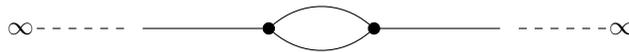

The {double bridge} $\mathcal B_2$ (Fig. \ref{ponte2}) is not semi--Eulerian, and the problem of
establishing the existence or the nonexistence of ground states becomes much more difficult than in $\mathcal B_3$.

\noindent
However, once again, as we shall see
$$
\E_{\mathcal B_2} (\mu)  \ = \ \E_\R (\mu)
$$
and  there is {\em no ground state}.

\noindent
We can then conclude that
on {bridge--graphs} the infimum is {never achieved}.


 \section{A key assumption} \label{XXXsec4}
From the examples of existence and of nonexistence given in the last section, one can single out
some observations:
\begin{enumerate}
\item It seems convenient to escape "intricated" zones,
e.g. vertices with high degree: at least, this is what happens for the star-graphs and for the bridges.

\item It is always possible to construct a soliton (possibly at
infinity, as the asymptotics of a sequence), so, in order to be a ground state,
a function must 
reach an energy level which is lower than the level of the soliton.

\item Therefore, in order to ensure existence 
of a ground state, the graph must exhibit
{\em structures}
able to trap functions that do better than the soliton.

\item
On the other hand, in order not to have minimizers, it seems
sufficient for the graph to be, in some sense, {\em more intricated than a line}.

\end{enumerate}
The last observation is  embodied in a topological assumption, that we call (H). We give
three alternative formulations of such assumption. The proof of the equivalence of the three formulations is not completely straightforward,
and will not be given here.

\begin{figure}[t]
\begin{center}
\begin{tikzpicture}
\node at (0,0) [nodo] (1) {};
\node at (-1.5,0) [infinito]  (2){$\infty$};
\node at (1,0) [nodo] (3) {};
\node at (0,1) [nodo] (4) {};
\node at (-1.5,1) [infinito] (5) {$\infty$};
\node at (2,0) [nodo] (6) {};
\node at (3,0) [nodo] (7) {};
\node at (2,1) [nodo] (8) {};
\node at (3,1) [nodo] (9) {};
\node at (4.5,0) [infinito] (10) {$\infty$};
\node at (5.5,0) [infinito] (11) {$\infty$};
\node at (4.5,1) [infinito] (12) {$\infty$};
\node at (1.5,0.85) {$\scriptstyle x$};
\node at (1.4,0.7) [nodino]  {};

\draw [-] (1) -- (2) ;
\draw [-] (1) -- (3);
\draw [-] (1) -- (4);
\draw [thick] (3) -- (4);
\draw [thick] (5) -- (4);
\draw [thick] (3) -- (6);
\draw [thick] (6) -- (7);
\draw [-] (6) to [out=-40,in=-140] (7);
\draw [thick] (3) to [out=10,in=-35] (1.4,0.7);
\draw [thick] (1.4,0.7) to [out=145,in=100] (3);
\draw [-] (6) to [out=40,in=140] (7);
\draw [-] (6) -- (8);
\draw [-] (6) to [out=130,in=-130] (8);
\draw [-] (7) -- (8);
\draw [-] (8) -- (9);
\draw [-] (7) -- (9);
\draw [-] (9) -- (12);
\draw [thick,red] (7) -- (10);
\draw [-] (7) to [out=40,in=140] (11);
\end{tikzpicture}
\end{center}
\caption{A point $x$ in the graph and a trail containing $x$ and two halflines.}
\label{trail}
\end{figure}
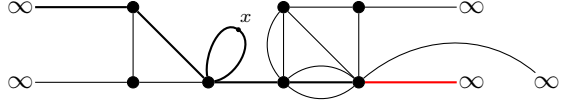

  The first formulation is based on the graph-theoretical notion of {\em trail} (Fig. \ref{trail}).
A  path made of adjacent edges, in which  every edge 
is run through exactly once, is called a {trail}. Notice that in a trail vertices can be
run through more than once.

\medskip 
  \noindent
{\bf Assumption (H), first formulation.}
{\em Every $x\in \G$ lies on a {\em trail} 
{that contains \em{two} halflines}}.

\begin{figure}[t]
\begin{center}
\begin{tikzpicture}
\node at (0,0) [nodo] (1) {};
\node at (-1.5,0) [infinito]  (2){$\infty$};
\node at (1,0) [nodo] (3) {};
\node at (0,1) [nodo] (4) {};
\node at (-1.5,1) [infinito] (5) {$\infty$};
\node at (2,0) [nodo] (6) {};
\node at (3,0) [nodo] (7) {};
\node at (2,1) [nodo] (8) {};
\node at (3,1) [nodo] (9) {};
\node at (4.5,0) [infinito] (10) {$\infty$};
\node at (5.5,0) [infinito] (11) {$\infty$};
\node at (4.5,1) [infinito] (12) {$\infty$};

\draw [-] (1) -- (2) ;
\draw [-] (1) -- (3);
\draw [-] (1) -- (4);
\draw [-] (3) -- (4);
\draw [-] (5) -- (4);
\draw [-] (6) -- (7);
\draw [-] (6) to [out=-40,in=-140] (7);
\draw [-] (3) to [out=10,in=-35] (1.4,0.7);
\draw [-] (1.4,0.7) to [out=145,in=100] (3);
\draw [-] (6) to [out=40,in=140] (7);
\draw [-] (6) -- (8);
\draw [-] (6) to [out=130,in=-130] (8);
\draw [-] (7) -- (8);
\draw [-] (8) -- (9);
\draw [-] (7) -- (9);
\draw [-] (9) -- (12);
\draw [-] (7) -- (10);
\draw [-] (7) to [out=40,in=140] (11);
\end{tikzpicture}
\end{center}
\caption{Once removed the only finite cut-edge of the graph, each connected component contains a halfline, and therefore
a vertex at infinity.}
\label{form2}
\end{figure}
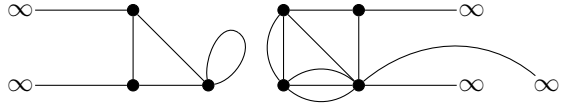

\begin{figure}[t]
\begin{center}
\begin{tikzpicture}
\node at (0,0) [nodo] (1) {};
\node at (-1.5,0) [infinito]  (2){$\infty$};
\node at (1,0) [nodo] (3) {};
\node at (0,1) [nodo] (4) {};
\node at (-1.5,1) [infinito] (5) {$\infty$};
\node at (2,0) [nodo] (6) {};
\node at (3,0) [nodo] (7) {};
\node at (2,1) [nodo] (8) {};
\node at (3,1) [nodo] (9) {};
\node at (4.5,0) [infinito] (10) {$\infty$};
\node at (5.5,0) [infinito] (11) {$\infty$};
\node at (4.5,1) [infinito] (12) {$\infty$};
\draw [-] (1) -- (3);
\draw [-] (1) -- (4);
\draw [-] (3) -- (4);
\draw [-] (5) -- (4);
\draw [-] (3) -- (6);
\draw [-] (6) -- (7);
\draw [-] (6) to [out=-40,in=-140] (7);
\draw [-] (3) to [out=10,in=-35] (1.4,0.7);
\draw [-] (1.4,0.7) to [out=145,in=100] (3);
\draw [-] (6) to [out=40,in=140] (7);
\draw [-] (6) -- (8);
\draw [-] (6) to [out=130,in=-130] (8);
\draw [-] (7) -- (8);
\draw [-] (8) -- (9);
\draw [-] (7) -- (9);
\draw [-] (9) -- (12);
\draw [-] (7) -- (10);
\draw [-] (7) to [out=40,in=140] (11);
\end{tikzpicture}
\end{center}
\caption{Assumption (H) in the second formulation applies to halflines too:
After removing a halfline, two connected components result: one contains some halflines, and therefore vertices at infinity; the
other is made of one vertex at infinity.}
\label{form2bis}
\end{figure}
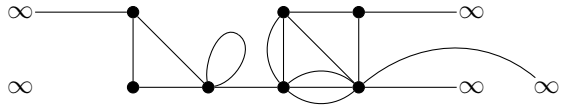

 \medskip \noindent
Assumption (H) can also be expressed as the absence of structure like ''bottlenecks'':
\medskip 
  
  \noindent
{\bf Assumption (H), second formulation.}
{\em After removing an arbitrary edge from $\G$,
every resulting connected component
contains a vertex at infinity} (Fig. \ref{form2}, \ref{form2bis}).



\medskip

  \noindent
The last formulation we give is more pictorial and considers the possibility
of covering the graph (vertex at infinity included) by cycles.

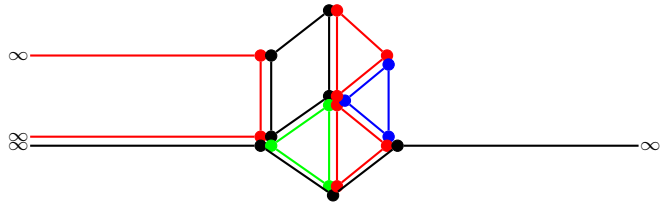
\begin{figure}[ht]
\begin{center}
\begin{tikzpicture}[xscale= 0.7,yscale=0.6]

\node at (-1.4,3) [nodo,red] (a){};
\node at (-1.2,3) [nodo] (b){};
\node at (-1.4,1.2) [nodo,red] (c){};
\node at (-1.2,1.2) [nodo] (d){};
\node at (-1.4,1) [nodo] (e){};
\node at (-1.2,1) [nodo,green] (f){};
\node at (-.1,4) [nodo] (g){};
\node at (-.1,2.1) [nodo] (h){};
\node at (-.1,1.9) [nodo, green] (i){};
\node at (-.1,.1) [nodo, green] (l){};
\node at (-.025,-.1) [nodo] (m){};
\node at (.05,.1) [nodo,red] (n){};
\node at (.05,1.9) [nodo,red] (o){};
\node at (.2,2) [nodo,blue] (p){};
\node at (.05,2.1) [nodo,red] (q){};
\node at (.05,4) [nodo,red] (r){};
\node at (1,3) [nodo,red] (s){};
\node at (1.03,2.8) [nodo,blue] (t){};
\node at (1.03,1.2) [nodo,blue] (u){};
\node at (1,1) [nodo,red] (v){};
\node at (1.2,1) [nodo] (w){};

\node at (-6,1.2) [infinito] (1) {$\scriptstyle\infty$};
\node at (-6,3) [infinito] (2) {$\scriptstyle\infty$};
\node at (-6,1) [infinito] (0) {$\scriptstyle\infty$};
\node at (6,1) [infinito] (3) {$\scriptstyle\infty$};

\draw[-,thick] (0)--(e);
\draw[-,red,thick] (1)--(c);
\draw[-,red,thick] (2)--(a);
\draw[-,thick] (3)--(w);

\draw[-,red,thick] (a)--(c);
\draw[-,thick] (b)--(d);
\draw[-,thick] (d)--(h);
\draw[-,thick] (h)--(g);
\draw[-,thick] (g)--(b);

\draw[-,green,thick] (f)--(l);
\draw[-,green,thick] (l)--(i);
\draw[-,green,thick] (i)--(f);

\draw[-,red,thick] (r)--(s);
\draw[-,red,thick] (s)--(q);
\draw[-,red,thick] (q)--(r);

\draw[-,blue,thick] (p)--(u);
\draw[-,blue,thick] (u)--(t);
\draw[-,blue,thick] (t)--(p);

\draw[-,red,thick] (v)--(o);
\draw[-,red,thick] (o)--(n);
\draw[-,red,thick] (n)--(v);

\draw[-,thick] (e)--(m);
\draw[-,thick] (m)--(w);

\end{tikzpicture}
\end{center}
\caption{A possible cycle covering.}
\label{cycle}
\end{figure}

\medskip 
  \noindent
{\bf Assumption (H), third formulation.}
  {\em After identifying all vertices at infinity,
the graph $\G$ admits a {\em cycle covering}} (Fig. \ref{cycle}).

\subsection{If (H) is violated}
Assumption (H) can be violated in several ways. Notice indeed that (H)
implies that $  {\G}$ has at least   {two} vertices at infinity, so it is violated by 
every graph having less than two halflines.
Furthermore, it is immediately seen that  {(H)} is violated not only by having less than two
halflines, but also by the presence of a {\em terminal edge} or {\em pendant} (Fig. \ref{pendente}):

\medskip

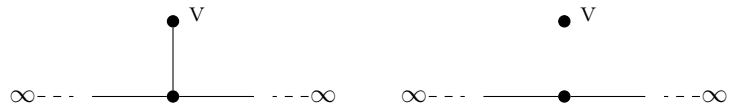
\begin{figure}
\begin{tikzpicture}[xscale=0.8]
\node at (-2.5,0) [infinito]  (1) {$\infty$};
\node at (0,0) [nodo] (2) {};
\node at (2.5,0) [infinito]  (3) {$\infty$};
\node at (1.5,0)  [minimum size=0pt] (5) {};
\node at (-1.5,0) [minimum size=0pt] (6) {};
\node at (0,1) [nodo] (7) {};
\node at (0.4,1.1) {$\vv$};

\draw[dashed] (1) -- (6);
\draw[dashed] (3) -- (5);
\draw [-] (2) -- (6) ;
\draw [-] (5) -- (2) ;
\draw [-] (2) -- (7) ;
\end{tikzpicture}
\qquad
\begin{tikzpicture}[xscale=0.8]
\node at (-2.5,0) [infinito]  (1) {$\infty$};
\node at (0,0) [nodo] (2) {};
\node at (2.5,0) [infinito]  (3) {$\infty$};
\node at (1.5,0)  [minimum size=0pt] (5) {};
\node at (-1.5,0) [minimum size=0pt] (6) {};
\node at (0,1) [nodo] (7) {};
\node at (0.4,1.1) {$\vv$};

\draw[dashed] (1) -- (6);
\draw[dashed] (3) -- (5);
\draw [-] (2) -- (6) ;
\draw [-] (5) -- (2) ;
\end{tikzpicture}
\caption{Left: line with a pendant. Right: after the removal of the pendant, there remains 
a connected compact component made of a vertex, violating (H) (see the second formulation).} 
\label{pendente}
\end{figure}

\begin{figure}
\begin{tikzpicture}[yscale=1,xscale=0.8]
\node at (3,0) [infinito]  (5) {$\scriptstyle\infty$}; 
\node at (5.5,0) [nodo] (6) {};
\node at (8,0) [infinito]  (7) {$\scriptstyle\infty$};
\node at (5.5,1) [nodo] (8) {};
\draw [-] (5) -- (6) ;
\draw [-] (6) -- (7) ;
\draw [-] (6) -- (8) ;
\draw(5.5,1.4) circle (0.4);
\node at (5.5,-.8)  [minimum size=0pt] (11){};
\end{tikzpicture}
\qquad
\begin{tikzpicture}[yscale=1,xscale=0.8]
\node at (3,0) [infinito]  (5) {$\scriptstyle\infty$}; 
\node at (5.5,0) [nodo] (6) {};
\node at (8,0) [infinito]  (7) {$\scriptstyle\infty$};
\node at (5.5,1) [nodo] (8) {};
\draw [-] (5) -- (6) ;
\draw [-] (6) -- (7) ;
\draw(5.5,1.4) circle (0.4);
\node at (5.5,-.8)  [minimum size=0pt] (11) {};
\end{tikzpicture}
\caption{Left: Signpost graph. Right: after the removal of the pendant, there remains 
a connected compact component made of a loop, violating (H) (see the second formulation).} 
\label{cartello}
\end{figure}
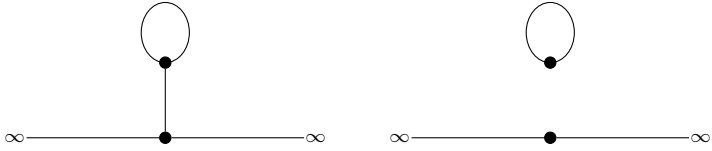

\noindent
 It is also possible to violate assumption (H) without having a pendant, like in the
 {\em signpost graph} (Fig. \ref{cartello}).

\section{A nonexistence result}
\begin{figure}[h] 
\begin{center}
\begin{tikzpicture}
\node at (-3,0) [infinito]  (1) {$\infty$};
\node at (0,0) [nodo] (2) {};
\node at (3,0) [infinito]  (3) {$\infty$};
\node at (2,0)  [minimum size=0pt] (5) {};
\node at (-2,0) [minimum size=0pt] (4) {};
\node at (0,1.2) [nodo] (6) {};
\node at (0,1.8) [nodo] (7) {};
\node at (0,0.6) [nodo] (8) {};

\draw[dashed] (4) -- (1);
\draw [-] (4) -- (2);
\draw [-] (5) -- (2) ;
\draw[dashed] (5) -- (3);
\draw(0,0.3) circle (0.3cm);
\draw[dashed]  (0,0.9) circle (0.3cm);
\draw[dashed]  (0,1.5) circle (0.3cm);
\draw (0,2.1) circle (0.3cm);
\end{tikzpicture}
\end{center}
\caption{A bubble tower.} \label{e3}
\end{figure}
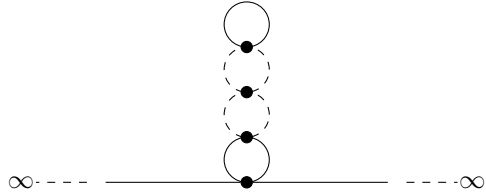

The first general result we give on graphs is negative:
\begin{theorem}[Nonexistence] \label{nonex}
Assume that $  \G$ satisfies assumption   {(H)}. Then, for any $\mu > 0$
$$
\EE_\G (\mu)  \ = \EE_\R (\mu)
$$
and the infium is   {never attained}, so that a ground state does not exist, except if $  \G$
is a ``bubble tower'' (Fig. \ref{e3}).

\end{theorem}

\noindent
The key idea of the theorem is strictly
related to {\em rearrangement theory}.

{  \subsection{Rearrangements}
The technique of rearrangements is nowadays classical in calculus of variations, and is 
widely used in order to show the existence of minimizers and possibly to establish some of their features, like for instance the 
symmetries. Its first extension to graphs is due
to L. Friedlander (\cite{friedlander}). 
In our proofs we use rearrangements to show the nonexistence of ground states. In what follows we give an intuitive and tutorial
summary of the results we will use, for readers non familiar with rearrangements. 

Given a nonnegative function ${u \in H^1_\mu (\G)}$, we aim at constructing
another nonnegative function ${v \in H^1 (\R^+)}$ s.t.
${E (v, \R^+) \leq E (u, \G)}.$
To this purpose, one can construct the so-called
{\em monotone rearrangement}. The idea behind it can be roughly summarized as
to
{\em cutting the graph of the function in vertical slices and  locating them on a halfline in
order of
decreasing height.}

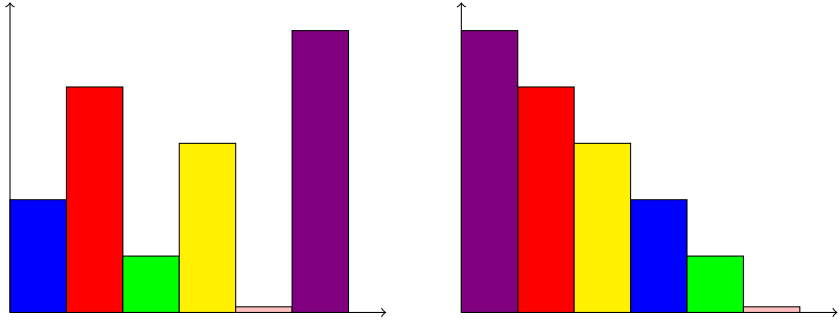
\begin{figure} 
\begin{tikzpicture}[yscale=0.75]
\draw [->](0,0) -- (5,0);
\draw [->](0,0)--(0,5.5);
\draw [fill=blue] (0,0) rectangle (0.75,2);
\draw [fill=red] (0.75,0) rectangle (1.5,4);
\draw [fill=green] (1.5,0) rectangle (2.25,1);
\draw [fill=yellow] (2.25,0) rectangle (3,3);
\draw [fill=pink](3,0) rectangle (3.75, 0.1);
\draw [fill=violet] (3.75,0) rectangle (4.5,5);

\draw [->](6,0) -- (11,0);
\draw [->](6,0)--(6,5.5);
\draw [fill=blue] (8.25,0) rectangle (9,2);
\draw [fill=red] (6.75,0) rectangle (7.5,4);
\draw [fill=green] (9,0) rectangle (9.75,1);
\draw [fill=yellow] (7.5,0) rectangle (8.25,3);
\draw [fill=pink](9.75,0) rectangle (10.5, 0.1);
\draw [fill=violet] (6,0) rectangle (6.75,5);

\end{tikzpicture}
\caption{Monotone rearrangement of a piecewise constant function: this example shows
pictorially the fact that integral norms are preserved while oscillations are reduced, so that the nonlinear
term in the energy is left untouched, while the kinetic term is dumped, so that the energy has diminished.}
\label{stogramma}
\end{figure}

As it appears from Fig.\ref{stogramma},
$L^p$-norms are preserved, while   {oscillations} are suppressed, so that one can
argue that, generalizing the procedure to regular functions, after rearranging a function {\em 
the kinetic
energy diminishes}.

Of course, one can give a more formal and general definition of monotone rearrangement:
let $(\Omega, {\mathcal F}, m)$ be a measure space, and consider a nonnegative function $f : \Omega \longrightarrow \R^+$.
One can define the
{\em distribution function} $\rho_f$ of the function $f$ as
$$ \rho_f (t) : = m ( \{ x \in \Omega, \ f (x) > t \} ). $$
Clearly, $\rho_f$ 
is defined $\R^+ \longrightarrow \R^+$ and is monotonically decreasing.
The monotone rearrangement $f^*$ of $f$ is a function defined on $\R^+$ with values in $\R^+$,
given by
$$ f^* (x) : = \inf \{ t \geq 0, \,  \rho_f (t) \leq x
\}.
$$
It is straightforward that, if $\rho_f$ is invertible, then
$$f^* = \rho^{-1}_f$$

Since $  {\rho_f = \rho_{f^*}}$, $f$ and $f^*$ have the same level sets, so that
$$ \| f \|_{L^p (X)} = \| f^* \|_{L^p (\R^+)} .$$

As an example, consider the measure space
$X = [- \pi, \pi]$, and the function $ {f (x) = |\sin x|}$. Then, one can directly compute $  {\rho_f (t) = 2 \pi - 4 {\rm arcsin} x}$
and  $f^* (x) = \cos (x/4)$ (Fig. \ref{sempio}, \ref{riarrangiata-ex}).




\begin{figure}
\begin{tikzpicture}[xscale=0.6,yscale = 0.5]
\begin{axis}[axis lines=middle,enlargelimits]
\addplot
[domain=-pi:pi,samples=100,smooth,thick,blue]
{abs(sin (360 * x / (2 * pi))};
\addplot
[domain=-pi:pi,samples=140,smooth,thick,red]
{0.707};
\addplot
[domain=-0.05:0.8,variable=\t,
 samples=40,smooth,thick,green]
(0.785,t);
\addplot
[domain=-0.05:0.8,variable=\t,
 samples=40,smooth,thick,green]
(-0.785,t);
\addplot
[domain=-0.05:0.8,variable=\t,
 samples=40,smooth,thick,green]
(2.255,t);
\addplot
[domain=-0.05:0.8,variable=\t,
 samples=40,smooth,thick,green]
(-2.255,t);
\addplot
[domain=-2.255:-0.785,variable=\t,
 samples=40,smooth,ultra thick,red]
(t,0);
\addplot
[domain=0.785:2.255,variable=\t,
 samples=40,smooth,ultra thick,red]
(t,0);
\end{axis}
\end{tikzpicture} \qquad
\ \qquad
\begin{tikzpicture}[xscale=0.6,yscale = 0.5]
\begin{axis}[axis lines=middle,enlargelimits]
\addplot
[domain=0:1,samples=40,smooth,thick,blue]
{2 * pi - 8 * pi * asin (x) / (360)};
\addplot
[domain=-0.05:4.2,variable=\t,
 samples=40,smooth,thick,red]
(0.707,t);
\addplot
[domain=-0.05:0.8,variable=\t,
 samples=40,smooth,thick,green]
(t,3.14);
\addplot
[domain=0:3.14,variable=\t,
 samples=40,smooth,ultra thick,red]
(0,t);
\end{axis}
\end{tikzpicture}
\caption{Left: $f(x) = | \sin x |$. Right: $\rho_f (t) = 2 \pi - 4 {\rm arcsin} x$.}
\label{sempio}
\end{figure}
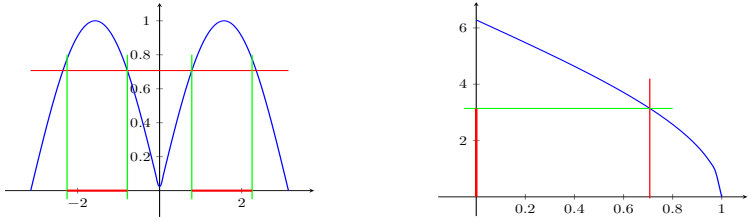


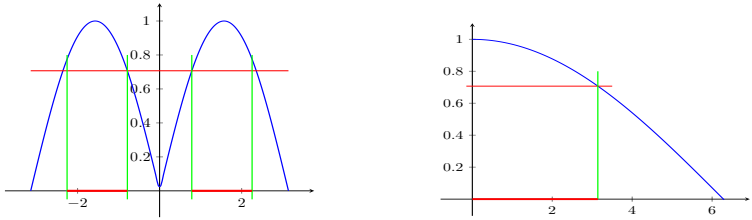
\begin{figure}
\begin{tikzpicture}[xscale=0.6,yscale = 0.5]
\begin{axis}[axis lines=middle,enlargelimits]
\addplot
[domain=-pi:pi,samples=100,smooth,thick,blue]
{abs(sin (360 * x / (2 * pi))};
\addplot
[domain=-pi:pi,samples=140,smooth,thick,red]
{0.707};
\addplot
[domain=-0.05:0.8,variable=\t,
 samples=40,smooth,thick,green]
(0.785,t);
\addplot
[domain=-0.05:0.8,variable=\t,
 samples=40,smooth,thick,green]
(-0.785,t);
\addplot
[domain=-0.05:0.8,variable=\t,
 samples=40,smooth,thick,green]
(2.255,t);
\addplot
[domain=-0.05:0.8,variable=\t,
 samples=40,smooth,thick,green]
(-2.255,t);
\addplot
[domain=-2.255:-0.785,variable=\t,
 samples=40,smooth,ultra thick,red]
(t,0);
\addplot
[domain=0.785:2.255,variable=\t,
 samples=40,smooth,ultra thick,red]
(t,0);
\end{axis}
\end{tikzpicture} \qquad
\ \qquad
\begin{tikzpicture}[xscale=0.6,yscale = 0.45]
\begin{axis}[axis lines=middle,enlargelimits]
\addplot
[domain=0:6.3,samples=40,smooth,thick,blue]
{cos(45* x/(pi))};
 \addplot
[domain=-0.15:3.5,variable=\t,
 samples=40,smooth,thick,red]
(t,0.707);
\addplot
[domain=0:0.8,variable=\t,
 samples=40,smooth,thick,green]
(3.14,t);
\addplot
[domain=0:3.14,variable=\t,
 samples=40,smooth,ultra thick,red]
(t,0);
\end{axis}
\end{tikzpicture}
\caption{Left:
${   f (x) = |\sin (x)|}$.
Right:
${   f^* (x) = \cos (x/4)}$
}
\label{riarrangiata-ex}
\end{figure}
 
\noindent
The main result that we borrow from rearrangement theory
consists in quantitatively estimating the decrease of the kinetic energy
induced by a rearrangement.

\noindent
To this aim, we first show
that the kinetic energy of the monotone rearrangement $u^* \in H^1 (\R^+)$
cannot exceed the energy of the  original function
$u \in H^1 (\G)$. We limit ourselves to the case of a function $u$ regular enough,
so that
$\rho_u$ is differentiable
and 
$\G$ can be partitioned in intervals $  {I_j = [a_j, b_j]}$ such that  $u$
is monotone in every $I_j$ (see Fig. \ref{monotonicity}). Then, it is easily seen that

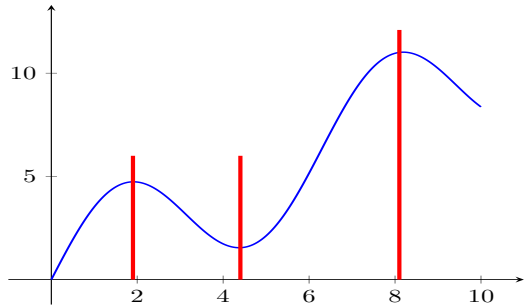
\begin{figure}
\begin{center}
\begin{tikzpicture}[xscale=1,yscale = 0.7]
\begin{axis}[axis lines=middle,enlargelimits]
\addplot
[domain=0:10,samples=100,smooth,thick,blue]
{x + 3 * sin (360 * x / (2 * pi))};
\addplot
[domain=0:6,variable=\t,
 samples=40,smooth,ultra thick,red]
(1.9,t);
\addplot
[domain=0:6,variable=\t,
 samples=40,smooth,ultra thick,red]
(4.4,t);
\addplot
[domain=0:12.1,variable=\t,
 samples=40,smooth,ultra thick,red]
(8.1,t);
\end{axis}
\end{tikzpicture}  
\end{center}
\caption{Partition of the domain of  a function  in intervals of monotonicity.}
\label{monotonicity}
\end{figure}

\begin{eqnarray*}
\rho (t+h) - \rho (t) \  & = & \ m  (\{ u(x) > t+h \}) - m (\{ u(x) > t \}) \\
& = & \sum_j \left( \ m ( \{x \in I_j, \ u(x) > t+h \}) - m ( \{x \in I_j, \ u(x) > t \} ) \right) \\
& \simeq & h \sum_j \ \f 1 {|u' (x_j (t))|},
\end{eqnarray*}
where $x_j(t)$ is the only point in $I_j$ where $u=t$. Then
$$ \rho'(t) =  \sum_{x, \, {\rm s.t.} \, u(x) = t} \ \f 1 {|u' (x)|}.$$

Thus, computing the kinetic energy  one finds
\begin{eqnarray*}
\int_\G | u' (x) |^2 \, dx & = & \sum_j \int_{a_j}^{b_j}
| u' (x) |^2 \, dx \
 = \ \sum_j \int_{\min_{[a_j, b_j]} u (x)}^{\max_{[a_j, b_j]} u(x)}
|u' (x_j(t))| \, dt \\
& = & \int_{0}^{\| u \|_\infty}
 \sum_{x, \, {\rm s.t.} \, u(x) = t} |u' (x)| \, dt
\end{eqnarray*}
where, in every interval $  {I_j}$, $  {t = u (x)}$.

\medskip
\noindent
Let $  {a_j > 0}$.  By Cauchy-Schwarz inequality,
\medskip
\begin{eqnarray*}
N & = & \sum_{j=1}^N 1 \ = \  \sum_{j=1}^N a_j^{1/2} a_j^{-1/2} \
 \leq \ \left( \sum_{j=1}^N a_j \right)^{1/2}
\left( \sum_{j=1}^N a_j^{-1} \right)^{1/2},
\end{eqnarray*}
so that, replacing $a_j$ with $|u'(x)|$,
\begin{equation*} \begin{split}
\int_\G | u' (x) |^2 \, dx \  & \geq \ \int_0^{\|u \|_\infty}
N^2 (t) \left( \sum_{x, \, {\rm s.t.} \, u(x) = t} \frac 1 {|u'(x)|} \right)^{-1}
\\ \ & = \ \int_0^{\|u \|_\infty} N^2 (t) \frac 1 {|\rho_u' (t) |},
\end{split} \end{equation*}
where, for every $t$ in the range of $u$, we defined the {\em number of preimages of $t$}
 $$N(t) : = \sharp u^{-1} (t).$$

\noindent    
Now, as $u^* = \rho^{-1}$, one gets
\begin{eqnarray*}
\int_\G | u' (x) |^2 \, dx \ & \geq & \ \int_0^{\|u \|_\infty}
\frac {N^2 (t)}  {|\rho_u' (t) |} \, dt \ \geq \  \int_0^{\|u \|_\infty}  \frac {dt} {|\rho_u' (t) |} \\
& = & \int_0^{\|u \|_\infty} |(u^*)' (x (t)) | \, dt \\
& = & \int_{\R^+} |(u^*)' (x ) |^2 \, dx
\end{eqnarray*}
where equality holds iff $  {N(t) = 1}$ for almost every $t$.

\noindent
We then proved the
{\em P\'olya-Szeg\H{o} inequality}:
\begin{equation} \nonumber 
\| (u^*)' \|_{L^2 (\R^+)} \, < \, \| u' \|_{L^2 (\G)}.
\end{equation}
Therefore, as
  the monotone rearrangement lowers the kinetic energy and preserves the nonlinear term, one has
$$ E (u, \G ) \ \geq \ E (u^*, \R^+).$$

\noindent
Notice that this implies that {the halfline is optimal} among non-compact graphs:
$$ \mathcal E_{\R^+} (\mu) \ \leq \mathcal E_\G (\mu), $$
for every graph $\G$ containing at least one halfline.

\subsection{Symmetric rearrangement}
To prove Theorem \ref{nonex} we need to introduce
the
{\em symmetric rearrangement} (Fig. \ref{symme}), defined as
$$ \widehat u (x) \ : = \ u^* (2 |x|), \qquad x \in \R$$

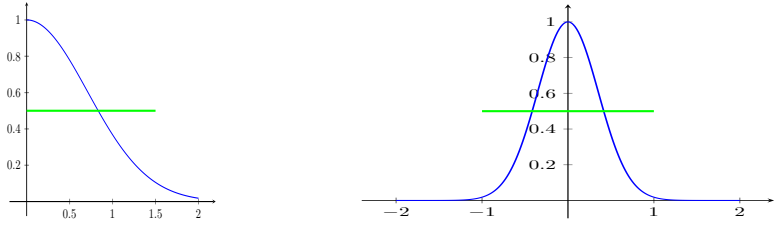
\begin{figure}
\begin{tikzpicture}[xscale=0.4,yscale=0.5]
\begin{axis}[axis lines=middle,enlargelimits]
\addplot
[domain=0:2,samples=100,smooth,thick,blue]
{e^(-x^2)};
\addplot
[domain=0:1.5,variable=\t,
 samples=40,smooth,ultra thick,green]
(t,0.5);
\end{axis}
\end{tikzpicture}
\qquad \ \, \ \qquad
\begin{tikzpicture}[xscale=0.8,yscale=0.5]
\begin{axis}[axis lines=middle,enlargelimits]
\addplot
[domain=-2:2,samples=100,smooth,thick,blue]
{e^(-4 * x^2)};
\addplot
[domain=-1:1,variable=\t,
 samples=40,smooth,ultra thick,green]
(t,0.5);
\end{axis}
\end{tikzpicture}
\caption{Left: The monotone rearrangement $u^*$. Right: The symmetric rearrangement $\widehat u$.}
\label{symme}
\end{figure}

\noindent
One immediately has that $\widehat u$ is even and
$$ \rho_{\widehat u} = \rho_{u^*} \ ( = \ \rho_u). $$

\noindent
By an elementary change of variable,
$$ \int_\R \widehat u^p \, dx \ = \ \int_{\R^+} (u^*)^p \, dx  \ \left( \ = \ 
\int_\G u^p \, dx \ \right)$$
and, analogously to the case of the monotone rarrangement, one finally has
\begin{eqnarray*} 
  \int_\R (\widehat u')^2 \, dx & = & 4 \int_{\R^+} [(u^*)']^2 \, dx \
= \ 4 \int_0^{\| u \|_\infty} \frac{dt}{|\rho'_u(t)|} \, dt \\
& \leq & \int_0^{\| u \|_\infty} \frac{N(t)^2}{|\rho'_u(t)|} \, dt \
= \ \int_\G (u')^2 \, dx
\end{eqnarray*}
provided that $  {N(t) \geq 2}$ for almost every $t$.

We finally proved the following
\begin{proposition}
Let $  {\G}$ be a connected non-compact metric graph, and $u$ 
be a nonnegative function in ${H^1(\G)}$.
Then, denoted
$$
N(t):=\#\{x\in \G\,:\,\, u(x)=t\}, \qquad t\in (0,\max u],
$$
the following inequality holds true:
$$
\int_{\R^+} |(u^*)'|^2\,dx\leq\int_{\G} |u'|^2\,dx,
$$
with strict inequality unless $  {N(t)=1}$ almost everywhere.

\noindent
Moreover, if
$
N(t)\geq 2$  almost everywhere, then
$$
\int_{\R} |(\widehat{u})'|^2\,dx \ \leq \ \int_{\G} |u'|^2\,dx,
$$
where equality implies that $  {N(t)=2}$ almost everywhere, and thus
$$ E (u, \G) \ \geq \ E (\widehat u, \R) \ \geq \ \mathcal E_\R (\mu)
\ = \ E (\phi_\mu, \R).
$$
\end{proposition}

%

%

%

We are now ready to prove the Theorem \ref{nonex}.
  
\medskip
\noindent
{\em Proof of Theorem \ref{nonex}.}
Let $  {u\in \Hmu(\G)}$, and let $  {x_0}$ be a   {global maximum point} for $   u$.

\noindent
Owing to Assumption (H) (see e.g. the first formulation), there exists
 a trail $  {\mathcal T}$ passing through $  { x_0}$ and containing two halflines.
Clearly, the restriction of $  {u}$ to $  {{\mathcal T}}$ belongs to $  {H^1({\mathcal T})}$ and $  {\max_{\mathcal T} u = \max_\G u}$.

\noindent
Furthermore, 
since $  {{\mathcal T}}$ connects   {two} vertices at infinity and $  {x_0 \in {\mathcal T}}$,
$$
\#\{x\in \G\,:\, u(x)=t\} \ge \#\{x\in{\mathcal T}\,:\, u(x)=t\}  { \ge 2}\quad \text{for \, a.e.}\; t.
$$
Due to the proposition on symmetric rearrangement, and to the
existence of runaway soliton sequences mimicking the soliton, the infimum can be
attained by a function $   u$, namely a ground state may exist, if and only if
\begin{enumerate}
\item Almost every point in $  {Ran \, u}$ has exactly two preimages.
\item $  {E(u, \G) = E (\phi_\mu, \R)}$
\end{enumerate}
Suppose that such a ground  state $   u$ exists, and call $  {x_0}$ a maximum point of $u$.
Then,
by assumption (H), there is a trail $\mathcal T$ passing through $x_0$
and containing two halflines. On this trail {\em every value in
$Ran \, u$ is attained twice}.

If there were other edges starting from (or arriving to) the trail,
then further counterimages would be created, and some interval
in $Ran \, u$ would be made of points with at least {\em three}
preimages, so that
$$ E (u, \G) \ > \ \mathcal E_\G (\mu) $$
contradicting the hypothesis of $u$ being a ground state.

Then, $  {\G = \mathcal T}$, i.e., $\G$ must be the real line (up to some possible identification of vertices),
and $u$ must be a soliton.

\noindent
The only identification of vertices that preserve the symmetry of the
soliton gives rise to the family of {\em tower of bubbles}. For an extended explanation of this point, see \cite{ast15-2}.

\hfill $\square$

\section{Ground states in the subcritical case}

In this section, given an exponent $p\in (2,6)$ and a mass $\mu>0$, we continue our study on the existence
of absolute minimizers (ground states) for the functional
\begin{equation} \nonumber
E (u, {\mathcal G}) 
\quad =\quad\frac12\int_\G |u'|^2\,dx - \frac1p \int_\G |u|^p\,dx,
\end{equation}
subject to the \emph{mass constraint}
\begin{equation} \nonumber
\int_{\G} |u|^2\,dx=\mu.
\end{equation}
Here $\G$ is an arbitrary noncompact metric graph (see Figure~\ref{how}), and
the range $(2,6)$ for the exponent $p$ is called the ``subcritical case'' (see Section \ref{XXXsecENRICO}
for the critical case where $p=6$).


Therefore, in trying to investigate ground states,
we shall be concerned with the case
where $\G$ does \emph{not} satisfy assumption (H).

A first result, regardless of ground states, is
that
the \emph{ground state energy level} is always \emph{intermediate} between
the half-soliton's (on the real halfline) and the soliton's (on real the line)
of the same mass $\mu$. More precisely, we have the following
\begin{theorem}[Level-pinching] \label{pinch}
For every non--compact graph $\G$,
$$
 E(\phi_{2\mu},\R^+)\quad \le\quad
 \inf_{v\in {H_\mu^1}(\G)} E(v,\G) \quad \leq\quad  E(\phi_\mu,\R)
$$
\end{theorem}

The first inequality is due to rearrangements: as explained in Sec. 5.1, given ${v\in H^1(\G)}$,
its decreasing rearrangement ${v^*}$ (over ${\R^+}$) has a lower
(possibly equal) energy.
In other words, no function $ v$ on $\G$ can ever beat the half-soliton on ${\R^+}$.

The {second inequality} (as explained in Section~\ref{XXXsec4}) is due
to the possibility of constructing {``quasi-solitons''}
escaping at $\infty$, along any half-line of $\G$ (since
$\G$ is noncompact, at least one of its edges must be unbounded, i.e. a half-line).
More precisely,
$\G$ contains arbitrarily
large intervals (in any half-line), and these intervals
can be used to support functions arbitrarily close to a soliton of mass $\mu$.

\bigskip

Of particular relevance is the case
where the \emph{second inequality is strict}.

\begin{theorem}[Existence of ground states] If ${\G}$ is non--compact and
$$
\inf_{v\in H^1_\mu(\G)} E(v,\G) {\quad<\quad} E(\phi_\mu,\R),
$$
then the infimum in \eqref{superE} is {attained}, i.e.  ${\G}$ supports a {ground state}.
\end{theorem}

Observe that the non--strict inequality {``$\leq$''} is always satisfied, due
to the level--pinching inequality (Thm. \ref{pinch}).

\medskip

The proof (see \cite{ast16-1}) is quite delicate, and is based on the
following \emph{dichotomy principle} for minimizing sequences
(relative to the infimum in \eqref{superE}). It turns out that, in general,
any minimizing sequence ${\{u_n\}}$ is \emph{either}
\begin{enumerate}
\item[(i)]   weakly convergent to zero, \emph{or}
\item [(ii)] strongly convergent to a ground state.
\end{enumerate}
But it can be proved that  {(i)} is (in this case)
\emph{incompatible} with the assumption of the theorem,
because ${u_n}$ would then ``escape to $\infty$'' along a halfline of $\G$,
approaching the shape of a soliton, and its energy level would then be \emph{equal}
to (and \emph{not} less than) the energy level of the soliton, in the limit.

\bigskip

The previous result is quite abstract, but it has the following consequence,
which is of quite practical use in the applications.

\begin{corollary} [Operative version of the existence theorem] If there exists
a competitor ${u\in H^1_\mu(\G)}$ such that
${E(u,\G) \, {\leq} \, E(\phi_\mu,\R)}$, then  ${\G}$ admits a {ground state}.
\end{corollary}

A sketch of the proof is as follows. Let  $ u$ be a competitor
satisfying the assumption of the theorem: if, by any chance, $ u$
is a {ground state}, then there is nothing to prove.
Otherwise $ u$ is not optimal, which amounts to
\[
\inf_{v\in H^1_\mu(\G)} E(v,\G) \,\,{<}\,\, E(u,\G) \leq E(\phi_\mu,\R),
\]
but in this case a ground state still exists (other than $u$) by the previous Theorem.

\bigskip

This corollary is quite useful in several concrete cases, where
one can try to obtain
estimates (on the ground state energy level) by  {graph surgery}:
starting from a soliton ${\phi_\mu}$ on {$\R$}, one can try to ``{fit it to $\G$}'',
without increasing its energy.
Whenever this can be done, the Theorem guarantees  that  $\G$ admits a {ground state}.

\medskip

A simple example where this can be done is the real line with a {pendant}, that is,
the graph in Fig.~\ref{figpend}.

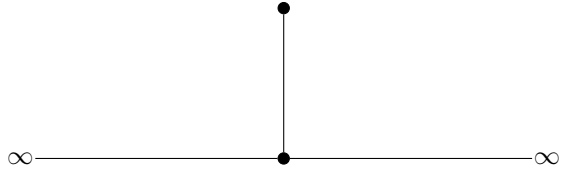
\begin{figure}
\begin{center}
\begin{tikzpicture}
\node at (-3.5,0) [infinito]  (1) {$\infty$};
\node at (0,0) [nodo] (2) {};
\node at (3.5,0) [infinito]  (3) {$\infty$};
\node at (0,2) [nodo] (7) {};
\draw [-] (2) -- (1) ;
\draw [-] (7) -- (2) ;
\draw [-] (2) -- (3) ;
\end{tikzpicture}
\end{center}
\caption{A line with one pendant (bounded edge) attached to it}
\label{figpend}
\end{figure}

\begin{theorem} Let ${\G}$ be the real line with a pendant of length $\ell$. Then
$$
\inf_{u\in H^1_\mu(\G)} E(u,\G) \,\,{<}\,\, E(\phi_\mu,\R),
$$
so that $\G$ admits a {ground state}.
\end{theorem}

The idea of the proof goes as follows. Due to the previous corollary,
it suffices to
construct a function ${u\in H^1_\mu(\G)}$ such
that ${E(u,\G)\,\,{<}\,\, E(\phi_\mu,\R)}$,
and this can be done
by \emph{graph surgery} combined with \emph{rearrangements}, as follows.
\begin{enumerate}
\item[(1)] Take the soliton ${\phi_\mu}$ centred
at zero and 	``cut it'' at a width ${\ell}$ (Fig. \ref{step1},\ref{step1bis}).
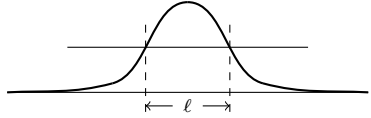
\begin{figure}
\begin{center}
\begin{tikzpicture}[xscale= 0.4,yscale=0.6]
\draw[-,thick] (0,2) to [out=0,in=170] (2.5,0.2);
\draw[-,thick] (0,2) to [out=180,in=10] (-2.5,0.2);
\draw[-,thick] (2.5,0.2) to [out=350,in=178] (6,0);
\draw[-,thick] (-2.5,0.2) to [out=190,in=2] (-6,0);
\draw[-] (-6,0) -- (6,0);
\draw (-4,1) -- (4,1);
\draw[dashed] (-1.4,1.5)--(-1.4,-0.5);
\draw[dashed] (1.4,1.5)--(1.4,-0.5);
\draw[<-] (-1.4, -0.3)--(-0.5,-0.3);
\draw[->] (0.5,-0.3)--(1.4, -0.3);
\node at (0,-0.3) {$\scriptstyle{\ell}$};
\end{tikzpicture}
\end{center}
\caption{First step: cut the head of the soliton}
\label{step1}
\end{figure}

\begin{figure}
\begin{center}
\begin{tikzpicture}[xscale= 0.58,yscale=0.8]
\draw[-,thick] (1,1.2) to [out=305,in=179] (5,0.02);
\draw[-,thick] (-1,1.2) to [out=235,in=2] (-5,0.02);
\draw[-|] (-5,0) -- (-1,0);
\draw[|-] (1,0) -- (5,0);
\node at (-1.2,-0.4) {{$\scriptscriptstyle -\ell/2$}};
\node at (1.2,-0.4) {{$\scriptscriptstyle\ell/2$}};
\draw[-,thick] (7,1) to [out=50,in=180] (8,1.8);
\draw[-,thick] (8,1.8) to [out=0,in=130] (9,1);

\draw[|-|] (7,0) -- (9,0);
\node at (6.8,-0.4) {{$\scriptscriptstyle -\ell/2$}};
\node at (9.2,-0.4) {{$\scriptscriptstyle\ell/2$}};
\node at (14,0) {};
\end{tikzpicture}
\end{center}
\caption{One is left with one head and two tails}
\label{step1bis}
\end{figure}
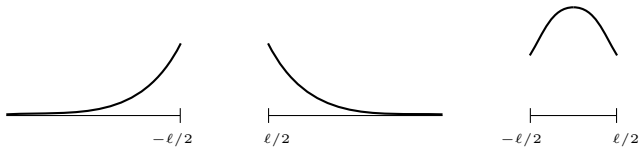

\item[(2)] Join the \emph{two} resulting soliton \emph{tails} at their maximum, and
  place them on the line in ${\G}$, with the maximum at
  the vertex (Fig. \ref{step2}).

\begin{figure}
\begin{center}
\begin{tikzpicture}[xscale= 0.6,yscale=0.8]
\draw[-,thick] (0,1.2) to [out=305,in=178] (4,0.02);
\draw[-] (0,0) -- (4,0);
\draw[-,thick] (0,1.2) to [out=235,in=2] (-4,0.02);
\draw[-|] (-4,0) -- (-0,0);
\node at (0,-0.4) {${\scriptstyle 0}$};
\end{tikzpicture}
\end{center}
\caption{Second step: glue the two tails together}
\label{step2}
\end{figure}
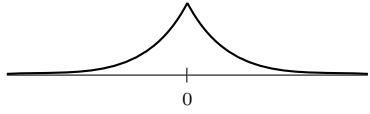

\item[(3)] {\emph{Rearrange}} the {head} of the soliton
to a
  {monotone function} on the interval {$[0, \ell]$} (Fig. \ref{step3}).

\begin{figure}
\begin{center}
\begin{tikzpicture}[xscale= 0.6,yscale=0.8]
\draw[|-|] (-2,0) -- (2,0);
\node at (-2,-0.4) {${\scriptstyle 0}$};
\node at (2,-0.4) {${\scriptstyle \ell}$};
\draw[-,thick] (-2,1) to [out=30,in=180] (2,2);
\draw[-,thick] (-9,1) to [out=50,in=180] (-7,2);
\draw[-,thick] (-7,2) to [out=0,in=130] (-5,1);
\draw[|-|] (-9,0) -- (-5,0);
\node at (-9.2,-0.4) {{$\scriptscriptstyle -\ell/2$}};
\node at (-5.2,-0.4) {{$\scriptscriptstyle\ell/2$}};
\node at (-3.5,1) {{$\Longrightarrow$}};
\end{tikzpicture}
\end{center}
\caption{Third step: rearrange monotonically the head of the soliton}
\label{step3}
\end{figure}
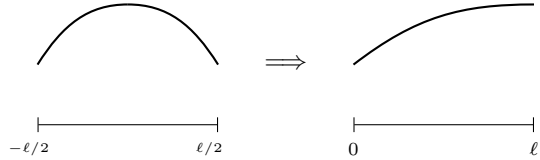

This monotone rearrangement {lowers the energy level}
of this portion of function.

\item[(4)] The function on the interval can be attached
to the function on the line, thus building a function on $\G$ (Fig. \ref{step4}):

\begin{figure}
\begin{center}
\begin{tikzpicture}[xscale= 0.6,yscale=0.8]
\node[nodo] at (2,1){};
\node at (-7.4,0) [infinito]   {$\infty$};
\node at (7.4,0) [infinito]   {$\infty$};

\draw[-,thick] (0,1.2) to [out=305,in=178] (4,0.02);
\draw[|-] (0,0) -- (7,0);
\draw[-,thick] (0,1.2) to [out=235,in=2] (-4,0.02);
\draw[-] (-7,0) -- (-0,0);

\draw[-,thick] (0,1.2) to [out=35,in=180] (2,2);
\node at (0,-0.4) {${\scriptstyle 0}$};
\node at (2.2,0.8){${\scriptstyle \ell}$};
\draw[-] (0,0) -- (2,1);

\draw[dotted] (2,1)--(2,2);
\draw[dotted] (0,0)--(0,1.2);

\end{tikzpicture}
\end{center}
\caption{Last step: mount the function on $\G$}
\label{step4}
\end{figure}
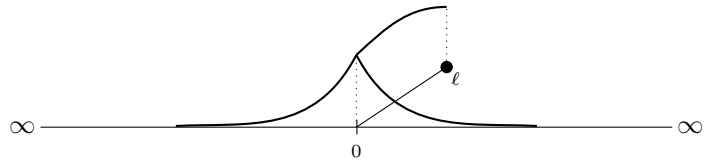

In this way, one produces a function ${u\in H^1_\mu(\G)}$ such that
$$
E(u,\G) < E(\phi_\mu,\R)
$$
(the strict inequality is due to the rearrangement
performed on the interval, ``from symmetric
to monotone'', in step (3)).
By the existence theorem, then, ${\G}$ admits a {ground state}.
\end{enumerate}

 We point out that we did \emph{not} construct the ground state,  but just
 a \emph{competitor} $u$, with an energy level
lower than the soliton's.

\bigskip

Other examples of graphs where the corollary
can be successfully applied are shown in Fig. \ref{towersign}, \ref{tad3}.

\bigskip

\begin{figure}
\begin{center}
\begin{tikzpicture}[scale= 1]

\node at (-3.5,0) [infinito]  (1) {$\scriptstyle\infty$};
\node at (-1,0) [nodo] (2) {};
\node at (1.5,0) [infinito]  (3) {$\scriptstyle\infty$};
\node at (-1,1) [nodo] (4) {};
\draw [-] (1) -- (2) ;
\draw [-] (2) -- (3) ;
\draw (-1,0.5) circle (0.5);
\draw (-1,1.3) circle (0.3);
\node at (3,0) [infinito]  (5) {$\scriptstyle\infty$}; 
\node at (5.5,0) [nodo] (6) {};
\node at (8,0) [infinito]  (7) {$\scriptstyle\infty$};
\node at (5.5,1) [nodo] (8) {};
\draw [-] (5) -- (6) ;
\draw [-] (6) -- (7) ;
\draw [-] (6) -- (8) ;
\draw(5.5,1.4) circle (0.4);
\end{tikzpicture}
\end{center}
\caption{Left: A line with a {tower of bubbles}. Right: A {signpost} graph. }
\label{towersign}
\end{figure}
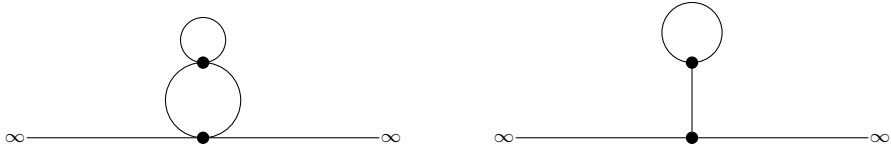

\bigskip

\begin{figure}
\begin{center}
\begin{tikzpicture}[scale= 0.9]
\node at (-2.5,0) [infinito]  (1) {$\scriptstyle\infty$};
\node at (1.8,0) [nodo] (2) {};
\draw [-] (1) -- (2) ;
\draw(2.2,0) circle (0.4);
\node at (4.3,0) [infinito]  (5) {$\scriptstyle\infty$}; 
\node at (8,0) [nodo] (6) {};
\node at (8.5,.5) [nodo]  (7) {};
\node at (9,-.6) [nodo] (8) {};
\node at (9.6,0) [nodo] (9) {};
\draw [-] (5) -- (6) ;
\draw [-] (6) -- (7) ;
\draw [-] (6) -- (8) ;
\draw [-] (6) -- (9) ;
\end{tikzpicture}
\end{center}
\caption{Left: A tadpole graph. Right: A 3-fork graph}
\label{tad3}
\end{figure}
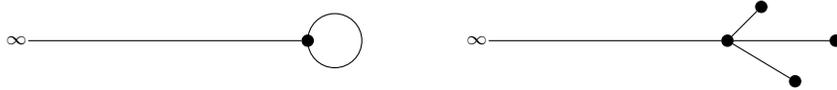

\bigskip

For each of these graphs, let us shortly see how one can build
a function ${u\in H^1_\mu(\G)}$ such
that ${E(u,\G)\leq  E(\phi_\mu,\R)}$, and thus prove the existence
of a ground state.

\bigskip

The first case, the so called ``bubble towers'', are graph of the kind
portrayed in Fig. \ref{bubbles} (as already seen in Theorem \ref{nonex}):

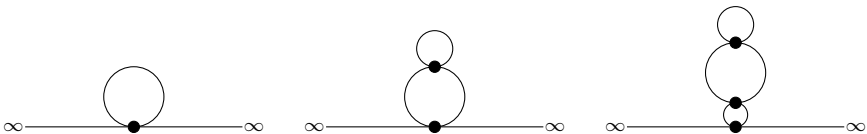
\begin{figure}
\begin{center}
\begin{tikzpicture} [xscale=.8,yscale=.8]

\node at(-7,0)  [infinito]  (1) {$\scriptstyle\infty$};
\node at(-3,0)  [infinito]  (3) {$\scriptstyle\infty$};
\node at(-5,0) [nodo] (2){};
\draw (-5,0.5) circle (0.5);
\draw [-] (1)--(2);
\draw [-] (2)--(3);

\node at(-2,0)  [infinito]  (4) {$\scriptstyle\infty$};
\node at(2,0)  [infinito]  (6) {$\scriptstyle\infty$};
\node at(0,0) [nodo] (5){};
\node at(0,1) [nodo] (10){};
\draw (0,0.5) circle (0.5);
\draw (0,1.3) circle (0.3);
\draw [-] (4)--(5);
\draw [-] (5)--(6);

\node at(3,0)  [infinito]  (7) {$\scriptstyle\infty$};
\node at(7,0)  [infinito]  (9) {$\scriptstyle\infty$};
\node at(5,0) [nodo] (8){};
\node at(5,.4) [nodo] (11){};
\node at(5,1.4) [nodo] (12){};
\draw (5,0.2) circle (0.2);
\draw (5,.9) circle (0.5);
\draw (5,1.7) circle (0.3);
\draw [-] (7)--(8);
\draw [-] (8)--(9);

\end{tikzpicture}
\end{center}
\caption{Some examples of bubble towers}
\label{bubbles}
\end{figure}

Each of them is obtained from ${\R}$, with the {identification} of some
pairs of {opposite points}:

\medskip

\centerline{${x_j}\,\,{\sim}\,\,{-x_j}$, $\quad{j=1,\ldots,n}$ \ \ ($ n$ bubbles)}

\medskip

The {symmetry} of these graphs enables them to {support a soliton} ${\phi_\mu}$,
exploting the {even symmetry} of the soliton:
\[
\phi_\mu(x_j)=\phi_\mu(-x_j),\qquad j=1,\ldots,n.
\]
As Fig. \ref{cutit} shows,
a soliton ${\phi_\mu}$  can indeed be
\emph{folded} and placed, \emph{isometrically}, on the line with two bubbles:

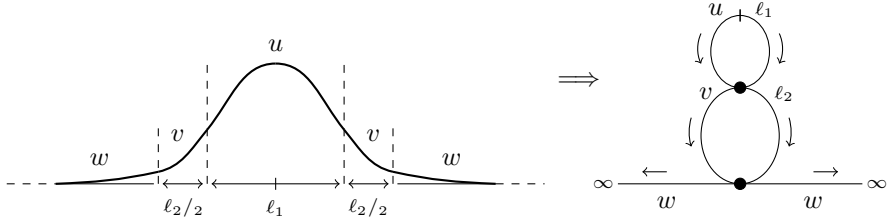
\begin{figure}
\begin{center}
\begin{tikzpicture}[xscale= 0.65,yscale=0.8]
\draw[-,thick] (-6,2) to [out=0,in=135] (-4.6,.9); 
\draw[-,thick] (-4.6,.9) to [out=315,in=170] (-3.6,0.2); 
\draw[-,thick] (-6,2) to [out=180,in=45] (-7.4,0.9); 
\draw[-,thick] (-7.4,.9) to [out=225,in=10] (-8.4,0.2); 
\draw[-,thick] (-3.6,0.2) to [out=350,in=178] (-1.5,0);
\draw[-,thick] (-8.4,0.2) to [out=190,in=2] (-10.5,0);
\node at (-8.4,0)  [infinito] (6) {\phantom{.}};
\draw [dashed] (-8.4,0) -- (-8.4,1);
\node at (-7.4,0)  [infinito] (7) {\phantom{.}};
\draw[dashed] (-7.4,0)--(-7.4,2);
\node at (-4.6,0)  [infinito] (8) {\phantom{.}};
\draw[dashed] (-4.6,0)--(-4.6,2);
\node at (-3.6,0)  [infinito] (9) {\phantom{.}};
\draw[dashed] (-3.6,0)--(-3.6,1);
\draw[<->] (6) -- (7);
\node at (0.2,1.7) {$\Longrightarrow$};
\draw[<->] (7) -- (8);
\draw[<->] (8) -- (9);
\draw[-] (9) -- (-1.5,0);
\draw[-] (6) -- (-10.5,0);
\draw[-] (-6,-0.1) -- (-6,0.1);
\draw[dashed] (-11.5,0) -- (-10.5,0);
\draw[dashed] (-1.5,0) -- (-.5,0);
\node at (-6,2.3)  {$u$};
\node at (-4,.8)  {${v}$};
\node at (-8,.8)  {${v}$};
\node at (-2.4,.4)  {${w}$};
\node at (-9.6,.4)  {${w}$};
\node at (-6,-0.4) {$\scriptstyle{\ell_1}$};
\node at (-4.1,-0.4) {$\scriptstyle{\ell_2/ _2}$};
\node at (-7.88,-0.4) {$\scriptstyle{\ell_2/ _2}$};
\draw[-] (1,0) -- (6,0);
\node at (.7,0) [infinito]  {$\scriptstyle\infty$};
\node at (6.3,0) [infinito]  {$\scriptstyle\infty$};
\node at (3.5,0) [nodo] (1) {};
\node at (3.5,1.6) [nodo] (2) {};
\draw (3.5,.8) circle (.8);
\draw(3.5,2.2) circle (.6);
\draw [->] (2.5,1.1) arc [radius=.7, start angle=160, end angle= 200];
\draw [->] (4.5,1.1) arc [radius=.7, start angle=20, end angle= -20];
\draw [->] (2.7,2.5) arc [radius=.5, start angle=150, end angle= 210];
\draw [->] (4.3,2.5) arc [radius=.5, start angle=30, end angle= -30];
\draw[->] (2,0.2) -- (1.5,0.2);
\draw[->] (5,0.2) -- (5.5,0.2);
\draw[-] (3.5,2.7) -- (3.5,2.9);
\node at (4,2.9) {$\scriptstyle{\ell_1}$};
\node at (4.4,1.5) {$\scriptstyle{\ell_2}$};
\node at (3,2.9) {$u$};
\node at (2.8,1.5) {${v}$};
\node at (2,-.3) {${w}$};
\node at (5,-.3) {${w}$};
\draw [-] (3.5,2.7) -- (3.5,2.9);  
\end{tikzpicture}
\end{center}
\caption{How to cut a soliton to fix it on a given bubble tower}
\label{cutit}
\end{figure}

Thus, in a sense, $\G$ ``supports'' a soliton ${\phi_\mu}$ and
this fact, combined with the level-pinching inequality, shows that
\[
\inf_{v\in H^1_\mu(\G)} E(v,\G)=E(\phi_\mu,\R).
\]
This {``folded soliton''} is therefore not just a competitor,
but precisely the \emph{ground state}.

In a similar way, one can see that any tower of bubbles
supports a (suitably folded) soliton, hence any tower of bubbles
has a ground state, and it is not difficult to show that the ground state
is unique, up to multiplication by a phase.

\bigskip

Also in the second example, the ``signpost graph'', there is  a ground state.
Indeed,
a soliton {$\phi_\mu$}, initially folded on a ``double bubble'', can be \emph{partially rearranged} and fitted to the signpost (see Fig. \ref{tosign})

\begin{figure}
\begin{center}
\begin{tikzpicture}[scale= 0.8]

\draw[-] (-6,0) -- (-1,0); 
\node at (-6.3,0) [infinito]  {$\scriptstyle\infty$};
\node at (-0.7,0) [infinito]  {$\scriptstyle\infty$};
\node at (-3.5,0) [nodo] (1) {};
\node at (-3.5,1.6) [nodo] (2) {};
\draw (-3.5,.8) circle (.8);
\draw(-3.5,2.2) circle (.6);
\node at (-4.3,2.6) {$u$};
\node at (-4.5,1.2) {$ v$};
\node at (-5,-.3) {$ w$};
\node at (-2,-.3) {$ w$};
\node at (0.2,1.7) {$\Longrightarrow$};

\node at (1,0) [infinito]  (5) {$\scriptstyle\infty$}; 
\node at (3.5,0) [nodo] (6) {};
\node at (6,0) [infinito]  (7) {$\scriptstyle\infty$};
\node at (3.5,2) [nodo] (8) {};
\draw [-] (5) -- (6) ;
\draw [-] (6) -- (7) ;
\draw [-] (6) -- (8) ;
\draw(3.5,2.6) circle (0.6);

\node at (2.7,3) {$u$};
\node at (3,1.1) {$ {v^*}$};
\node at (4.9,-.3) {$ w$};
\node at (2.1,-.3) {$ w$};

\end{tikzpicture}
\end{center}
\caption{From a bubble tower to a signpost}
\label{tosign}
\end{figure}
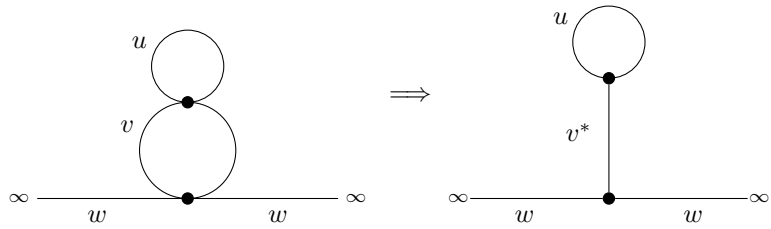
In the above picture, ${v^*}$ denotes the monotone rearrangement of $ v$
(from the circle to an interval of the same length, that is, ``from symmetric
to monotone'').
The {loss of preimages} in passing from $ v$ (regarded as an {even}
function) to
${v^*}$ (regarded as a {decreasing} function) makes the
energy decrease and go \emph{below} ${E(\phi_\mu,\R)}$.
As before,  we did {not} build a ground state, just a good {competitor}.

\bigskip

Also in the case of the ``tadpole graph''
 we can \emph{partially rearrange} the competitor alredy built
on the double bubble (see Fig. \ref{totadpole})

\begin{figure}
\begin{center}
\begin{tikzpicture}[scale= 0.9]

\draw[-] (-6,0) -- (-1,0); 
\node at (-6.3,0) [infinito]  {$\scriptstyle\infty$};
\node at (-0.7,0) [infinito]  {$\scriptstyle\infty$};
\node at (-3.5,0) [nodo] (1) {};
\node at (-3.5,1.6) [nodo] (2) {};
\draw (-3.5,.8) circle (.8);
\draw(-3.5,2.2) circle (.6);
\node at (-4.3,2.6) {$u$};
\node at (-4.5,1.2) {$ v$};
\node at (-5,-.3) {$ w$};
\node at (-2,-.3) {$ w$};

\node at (0.2,1.7) {$\Longrightarrow$};

\node at (3.5,0) [nodo] (6) {};   
\node at (3.5,-3) [infinito]  (7) {$\scriptstyle\infty$};
\node at (3.5,2) [nodo] (8) {};
\draw [-] (6) -- (8) ;
\draw (3.5,2.6) circle (0.6);
\draw [-] (7) -- (6);

\node at (2.7,3) {$u$};
\node at (3,1.1) {$ {v^*}$};
\node at (3,-1.3) {${ w^*}$};

\end{tikzpicture}
\end{center}
\caption{From bubble tower to tadpole}
\label{totadpole}
\end{figure}
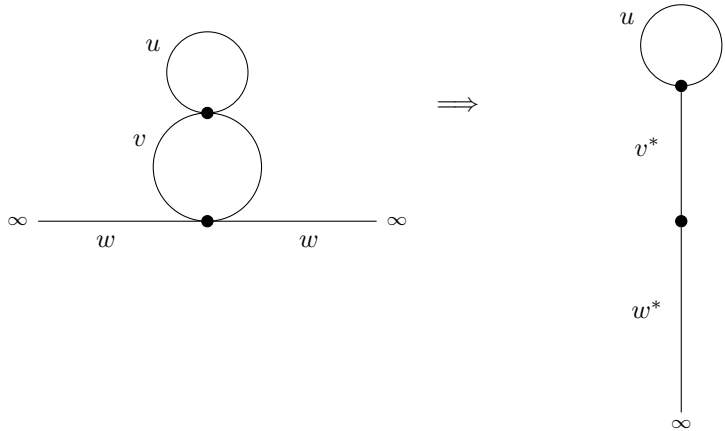

Here, in addition to the rarrangement of $v$, we also rearrange $w$
(from the real line) to $w^*$ (to the half-line), which further decreases energy.

\bigskip

Finally, a similar procedure applies to
the case where the graph is a ``3--fork''. Starting from the
competitor on a double bubble, we first ``open up'' the
two circles corresponding to the two bubbles,
and  we rearrange $w$ to $w^*$ (on the half-line),
as illustrated in Fig. \ref{openarc}

\medskip

\begin{figure}
\begin{center}
\begin{tikzpicture}[scale= 0.8]

\draw[-] (-7,0) -- (-2,0); 
\node at (-7.3,0) [infinito]  {$\scriptstyle\infty$};
\node at (-1.7,0) [infinito]  {$\scriptstyle\infty$};
\node at (-4.5,0) [nodo] (1) {};
\node at (-4.5,1.6) [nodo] (2) {};
\draw (-4.5,.8) circle (.8);
\draw(-4.5,2.2) circle (.6);
\node at (-5.3,2.6) {$u$};
\node at (-5.5,1.2) {$ v$};
\node at (-6,-.3) {$ w$};
\node at (-3,-.3) {$ w$};

\node at (-.8,1.7) {$\Longrightarrow$};

\node at (4.6,0) [nodo] {}; 
\node at (5.2,.6) [nodo] {};
\node at (5.6,.4) [nodo] {};
\node at (3,0) [nodo] {};
\node at (3.1,-.4)[nodo]{};
\node at (2.95,0) [infinito] (0) {};
\node at (0,0) [infinito] (1) {$\scriptstyle\infty$};
\draw [-] (0) -- (1);
\draw [-] (3,0) arc [radius=.8, start angle=180, end angle= -150];
\draw [-] (4.6,0) arc [radius=.6, start angle=180, end angle=90];
\draw [-] (4.6,0) arc [radius=.6, start angle=-180, end angle= 45];
\node at (1.7,-.3) {${w^*}$};
\node at (3.5,1) {$ v$};
\node at (5.1,.9) {$u$};
\node at (5.6,-.7) {$u$};
\node at (4,-1) {$ v$};
\end{tikzpicture}
\end{center}
\caption{Open the arcs in the bubble tower}
\label{openarc}
\end{figure}
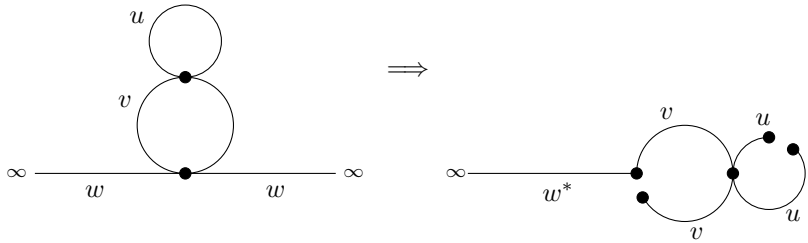

\medskip

Now, the arc of circle with the free endpoint from the lower bubble, and
the two arcs of circle from the upper bubble, can be seen as
the three bounded edges forming the fork (see Fig. \ref{pov}). Of course, by a proper
choice of the size of the bubbles and the cut-points, a 3-fork
with edges of any size (not necessarily equal) can be handled.

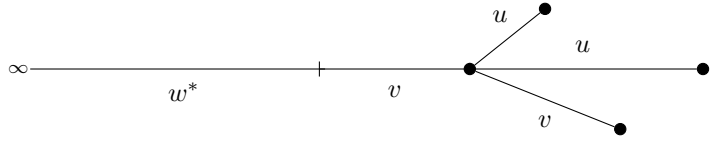
\begin{figure}
\begin{center}
\begin{tikzpicture}[scale= 1]

\node at (-6,-4) [infinito]  (1) {$\scriptstyle\infty$};
\node at (-3.8,-4.3) {${w^*}$};
\coordinate (2bis) at (0,-4);
\draw [-] (-2,-4.1) -- (-2,-3.9);

\node at (-1,-4.3) {$ v$};

\node at (-1.8,-4)  [minimum size=0pt] (10) {};

\draw[-] (1)--(10);
\draw[-] (-2,-4)--(0,-4);

\node at (0,-4) [nodo] (2) {}; 
\draw [-] (2bis) -- (2);
\node at (3.1,-4) [nodo]  (3) {};
\node at (-3.5,-4) [minimum size=0pt] (6) {};
\node at (1,-3.2) [nodo] (7) {};
\node at (2,-4.8) [nodo] (8) {};
\draw [-] (3) -- (2) ;
\draw [-] (2) -- (7) ;
\draw [-] (2) -- (8) ;
\node at (.4,-3.3) {$u$};
\node at (1.5,-3.7) {$u$};
\node at (1,-4.7) {$ v$};

\end{tikzpicture}
\end{center}
\caption{Reconstruct the $3$-fork}
\label{pov}
\end{figure}

  In the examples we have seen, the \emph{topology} of $\G$ was
enough to guarantee a ground state,
  while the {metric} of $\G$ (i.e. the {lengths} of its edges)
  was irrelevant.

  However, in general, things are more complicated, and also
  the \emph{metric} of $\G$ may {play a role}.

\medskip

  We will consider two examples where this is the case:

  \begin{itemize}
    \item graphs with just {one half-line};

    \item graphs where \emph{phase transitions} occur, from existence to
nonexistence of ground states,
if we {vary} the length of just \emph{one} edge.
  \end{itemize}

Let $\G$ be a graph with just \emph{one half-line} (Fig. \ref{just}).

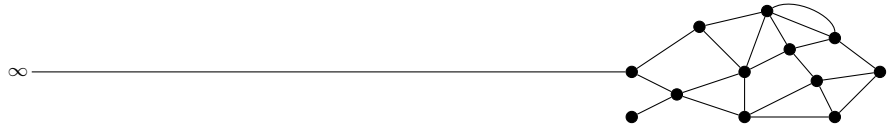
\begin{figure}
\begin{center}
\begin{tikzpicture}[xscale= 0.6,yscale=0.3]
\node at (-.5,2) [nodo] (02) {};
\node at (2,2) [nodo] (22) {};
\node at (2.5,4.7) [nodo] (24) {};
\node at (3.6,1.6) [nodo] (42) {};
\node at (3,3) [nodo] (33) {};
\node at (5,2) [nodo] (52) {};
\node at (4,3.5) [nodo] (43) {};
\node at (1,4) [nodo] (04) {};
\node at (.5,1) [nodo] (11) {};
\node at (-.5,0) [nodo] (11b) {};
\node at (2,0) [nodo] (20) {};
\node at (4,0) [nodo] (40) {};
\node at (-14,2) [minimum size=0pt] (meno) {};
\node at (-14.1,2) [infinito]  (infmeno) {$\scriptstyle\infty$};
\draw[-] (11)--(11b);
\draw[-] (02)--(04);
\draw[-] (04)--(24);
\draw[-] (04)--(22);
\draw[-] (24)--(22);
\draw[-] (02)--(11);
\draw[-] (11)--(22);
\draw[-] (11)--(20);
\draw[-] (20)--(22);
\draw[-] (22)--(33);
\draw[-] (24)--(33);
\draw[-] (24)--(43);
\draw[-] (24) to [out=45,in=90] (43);
\draw[-] (33)--(43);
\draw[-] (43)--(52);
\draw[-] (33)--(42);
\draw[-] (20)--(42);
\draw[-] (20)--(40);
\draw[-] (40)--(42);
\draw[-] (42)--(52);
\draw[-] (40)--(52);
\draw[-] (02)--(meno);
\end{tikzpicture}
\end{center}
\caption{A graph with one halfline}
\label{just}
\end{figure}

The main question is, of course, whether $\G$ admits
a ground state for every value of the prescribed mass {$\mu$}.
This is nontrivial, since one can check that
such a graph (due to the presence of \emph{just one} half-line)
does not satisfy assumption (H), and hence the existence of ground states
cannot be a priori ruled out.

As we have seen, several examples of graphs with just one half-line (tadpole, $2$--fork, $3$--fork) 
indeed admit a ground state for every  {$\mu$}.

\medskip

However, this is not true in general, and \emph{counterexamples} can be constructed.

Let {${\mathcal K}$} be {any} compact graph, and let {$\G$} be the graph obtained
by attaching \emph{one half-line} to {${\mathcal K}$} (Fig. \ref{just2}):
\begin{figure}
\begin{center}
\begin{tikzpicture}[xscale= 0.4,yscale=0.2]
\node at (-.5,0) [nodo] (11b) {};
\node at (-.5,2) [nodo] (02) {};
\node at (2,2) [nodo] (22) {};
\node at (2.5,4.7) [nodo] (24) {};
\node at (3.6,1.6) [nodo] (42) {};
\node at (3,3) [nodo] (33) {};
\node at (5,2) [nodo] (52) {};
\node at (4,3.5) [nodo] (43) {};
\node at (1,4) [nodo] (04) {};
\node at (.5,1) [nodo] (11) {};
\draw[-] (11)--(11b);
\node at (2,0) [nodo] (20) {};
\node at (4,0) [nodo] (40) {};
\node at (-14,2) [minimum size=0pt] (meno) {};
\node at (5.5,4) [minimum size=0pt] {{${\mathcal K}$}};
\node at (-14.1,2) [infinito]  (infmeno) {$\scriptstyle\infty$};
\draw[-] (24) to [out=45,in=90] (43);
\draw[-] (02)--(04);
\draw[-] (04)--(24);
\draw[-] (04)--(22);
\draw[-] (24)--(22);
\draw[-] (02)--(11);
\draw[-] (11)--(22);
\draw[-] (11)--(20);
\draw[-] (20)--(22);
\draw[-] (22)--(33);
\draw[-] (24)--(33);
\draw[-] (24)--(43);
\draw[-] (33)--(43);
\draw[-] (43)--(52);
\draw[-] (33)--(42);
\draw[-] (20)--(42);
\draw[-] (20)--(40);
\draw[-] (40)--(42);
\draw[-] (42)--(52);
\draw[-] (40)--(52);
\draw[-] (02)--(meno);
\end{tikzpicture}
\end{center}
\caption{A graph obtained by attaching a half-line to a compact graph $\mathcal K$.}
\label{just2}
\end{figure}
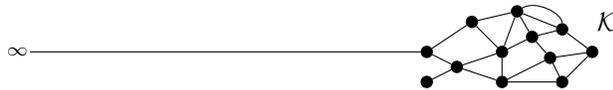
\begin{theorem}
There exists ${\varepsilon>0}$ such that if
\[
\mu^\beta \mathop{\rm diam} ({\mathcal K}) + \frac1{\mu^\beta \mathop{\rm length}({\mathcal K})}  <  \varepsilon, \qquad \beta = \frac{p-2}{6-p},
\]
then ${\G}$ has {no ground state} with mass ${\mu}$.
\end{theorem}

The proof is quite involved and requires several sharp estimates, the interested reader is referred to 
\cite{ast16-1}.
The main
idea, however, is simple: a \emph{small diameter} of the
compact core ${\mathcal K}$, combined with a \emph{long total length}, rules
 out ground
states, because 
\emph{tangled
compact parts} are not energetically convenient. Any competitor $u$,
due to
the structure of ${\mathcal K}$, 
is indeed forced to either \emph{oscillate} (and thus have many preimages) or,
on the contrary, to be \emph{almost constant}, and neither behaviour is 
energetically convenient if $u$ is compared to a soliton of the same mass.

\medskip

A concrete case where this result applies is when
$\G$ has the shape of an
\emph{$n$--fork}, namely $ n$ terminal edges
of length $\ell$, attached to a half-line (Fig. \ref{nfork}).

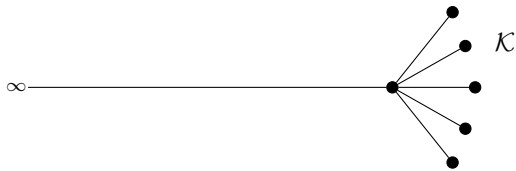
\begin{figure}
\begin{center}
\begin{tikzpicture}
\node at (-5,0) [infinito]  (1) {$\scriptstyle\infty$};
\node at (-4,0) [minimum size=0pt] (2) {};
\node at (1.5,.6) [minimum size=0pt]  {${{\mathcal K}}$};
\node at (0,0) [nodo] (3) {};
\node at (1.1,0) [nodo]  (4) {};
\node at (0.8,1) [nodo]  (5) {};
\node at (0.97,0.55) [nodo]  (6) {};
\node at (0.8,-1) [nodo]  (7) {};
\node at (0.97,-0.55) [nodo]  (8) {};
\draw [-] (1) -- (3) ;
\draw [-] (3) -- (4) ;
\draw [-] (3) -- (5) ;
\draw [-] (3) -- (6) ;
\draw [-] (3) -- (7) ;
\draw [-] (3) -- (8) ;
\end{tikzpicture}
\end{center}
\caption{The $n$-fork graph}
\label{nfork}
\end{figure}

In this case, we clearly have ${\mathop{\rm diam} (\mathcal K) = 2\, \ell}$\ \ and\ \
  ${\mathop{\rm length} (\mathcal K) = n \,\ell}$.
If we fix  the value of the mass {$ \mu>0$}, and
then we take {$\ell$} \emph{small enough} (depending on $\mu$)
 and {$n$} \emph{large enough} (depending on $\mu$ and $\ell$),
then  the theorem applies, and the resulting $\G$ has no ground state.

Explicit computations show that, at least when $p=4$, 
any  ${n\ge 5}$ is \emph{sufficient} for the counterexample, while
on the other hand
${n >3}$ is \emph{necessary}, because we know that any $3$--fork has a ground state.

Finally, it is not known whether one can build the counterexample
with ${n=4}$, that is,
it is not known whether a $4$-fork always has a ground state.

\bigskip

Now we discuss an example of a metric
graph $\G$ such that varying the length of just one edge
(without affecting the topology of $\G$) may lead from existence to
nonexistence of a ground state.

Let $\G_\ell$ consist of three half-lines and one \emph{terminal edge}  of length $\ell$,
all emanating from a common vertex (Fig. \ref{emanating})
\begin{figure}
\begin{center}
\begin{tikzpicture}
\node at (-3,0) [infinito]  (1) {$\scriptstyle\infty$};
\node at (0,0) [nodo] (2) {};
\node at (3,0) [infinito]  (3) {$\scriptstyle\infty$};
\node at (0,1)  (4) [nodo] {};
\node at (0,-1.5) [infinito]  (8) {$\scriptstyle\infty$};
\node at (.2,.6) {$\scriptstyle\ell$};
\draw [-] (1) -- (2) ;
\draw [-] (2) -- (3) ;
\draw [-] (2) -- (4) ;
\draw [-] (2) -- (8) ;
\end{tikzpicture}
\end{center}
\caption{A graph made of three half-lines and a terminal edge}
\label{emanating}
\end{figure}
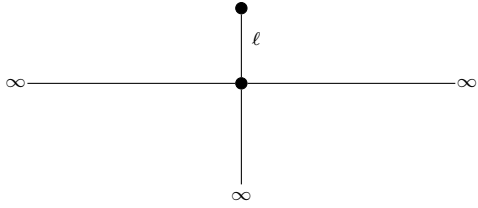

Clearly, as long as $\ell>0$, the {topology} of {$\G_\ell$} is \emph{independent} of the length
{$\ell$}. Nevertheless, we have the following

\begin{theorem}[phase transition] There exists a critical length
${\ell^*>0}$ such that:

\medskip

\centerline{$
\text{{$\G_\ell$} has a ground state}
\quad{\iff}\quad  {\ell\ge \ell^*}.
$}

\end{theorem}

The idea of the proof is that, once the mass $\mu$ has been fixed,
if $\ell$ is long enough then $\G_\ell$ has a ground state (this is true in general,
as soon as a graph $\G$ has a long enough terminal edge), which resembles a half-soliton
with the head at the tip of the bounded edge. On the other hand, if $\G_\ell$
had a ground state for every
$\ell>0$, then by a compactness argument also
the 3-star graph $\G_0$ would inherit a ground state, which is a contradiction
since $\G_0$ is known to have no ground state. This shows that at least one transition
(from existence to nonexistence of a ground state) 
must occur, as $\ell$ is decreased: then, the fact
that \emph{exactly one} transition occurs requires a more careful analysis, based on
a monotonicity argument.

This example shows that, in general, the \emph{topology} of $\G$ 
\emph{is not enough, alone,} to establish whether $\G$ has a ground state
of a given mass, and also the \emph{metric properties}
of $\G$ (together with its topology) should be considered.


\newcommand{\HH}{\mathcal H}
\newcommand{\Gp}{\mathcal G'}
\newcommand{\K}{\mathcal K}
\newcommand{\C}{\mathcal C}
\newcommand{\eps}{{\varepsilon}}

\newcommand\dx{{\,dx}}

\newcommand\mur{\mu_\R}
\newcommand\mup{\mu_{\R^+}}
\newcommand\mug{\mu_\G}


\renewcommand{\thefootnote}{\fnsymbol{footnote}}


\renewcommand{\thefootnote}{\arabic{footnote}}
\setcounter{footnote}{0}
%

\section{The critical case: $p=6$} \label{XXXsecENRICO}

In this section we describe some results  concerning  the existence of ground states for the {\em critical} NLS energy functional

\[
E(u,\G)= \frac12\int_\G |u'|^2\,dx - \frac16 \int_\G |u|^6\,dx
\]
on the space $H^1_\mu (\G)$, where $\G$ is a noncompact metric graph (see Fig. 1). The content of this section
refers to \cite{ast16-2}.

\medskip

\noindent
According to \eqref{el},
the solutions to \eqref{superE} are solutions of the $L^2$--critical stationary NLS equation
\[
u'' + u^5 = \omega u \qquad\text{ on } \; \G,
\]
with Kirchhoff conditions \eqref{kirch} at each vertex $\vv$ of the graph.


This problem is much more delicate than the subcritical one, where the exponent in the nonlinearity
lies in the interval  $(2,6)$.
One of the reasons is that, as discussed in Sec. 2.1 
under the formal mass-preserving transformation 
\[
u(x) \mapsto u_\lambda(x) = \sqrt\lambda u(\lambda x),
\]
the kinetic and the potential terms in $E$ scale in the same way:
\begin{equation} \label{homogeneity}
E(u_\lambda, \lambda^{-1}\G) = \lambda^2 E(u,\G),
\end{equation}
which is typical of problems with {\em serious} loss of compactness.

In the critical case  the problem depends very strongly on the mass  $\mu$ and on the 
ground state energy function
$
\EE_\G(\mu) 
$ defined in \eqref{superE},
which will play a central role in all of our results.
\medskip

\begin{example} The real line  $(\G = \R)$. The situation is very different from
the one encountered in Sec. 3.1 for the subcritical case.

\noindent Indeed,
it is known that there exists a number $\mu_\R>0$, the {\em critical mass}, such that
\[
\EE_\R(\mu) = \begin{cases} 0 & \text{ if } \mu \le \mu_\R \\  -\infty & \text{ if } \mu > \mu_\R \end{cases}\quad\qquad \left(\mur = \pi\sqrt 3 /2\right).
\]
Moreover $\EE_\R(\mu)$ is attained {\em if and only if} $\mu = \mu_\R$. 
The {\em ground states}, called solitons, form a quite large family: up to phase and translations, they can be written as
\[
\phi_\lambda(x) =  \sqrt\lambda \phi(\lambda x), \qquad \lambda>0,
\]
where $\phi (x) =  \hbox{\rm{sech}}^{1/2}(\frac2{\sqrt 3} x)$.
\end{example}

\begin{example} The half-line $(\G = \R^+)$.
Again, there exists a number $\mup = \mur/2$, such that
\[
\EE_{\R^+}(\mu) = \begin{cases} 0 & \text{ if } \mu \le \mup \\  -\infty & \text{ if } \mu > \mup \end{cases}\quad\qquad \left(\mup = \pi\sqrt 3 /4\right).
\]
Moreover $\EE_\R(\mu)$ is attained {\em if and only if} $\mu = \mup$. 
The ground states (half-solitons) are the restrictions to $\R^+$ of the family $\phi_\lambda$.





\end{example}

Thus on the standard domains $\R$ and $\R^+$ the minimization process \eqref{superE} is {\em extremely unstable},
with solutions existing for a {\em single value} of the mass.

This behavior is due to the same homogeneity of the kinetic and potential terms under mass-preserving scalings
and the invariance of $\R$ and $\R^+$ under dilations.

On a generic noncompact graph $\G$ however, the problem can be highly nontrivial and entirely new phenomena
may arise, depending on the topology of the graph.

Here we describe these new phenomena, essentially by classifying all graphs from the point of view
of existence of ground states.

\subsection{The critical mass}  

The first thing to do is to understand the appearance of the critical mass $\mur$ (or $\mup$) in
the problems on classical domains and to identify the same notion for general graphs. This is carried out
by analyzing the Gagliardo-Nirenberg inequality, a fundamental tool in all the existence proofs.

The Gagliardo--Nirenberg inequality on $\R$ reads
\[
\|u\|_6^6 \le C\|u\|_2^4\cdot \|u'\|_2^2 \quad \qquad \forall u\in H^1(\R).
\]
The {\em best constant} (the smallest $C$) is
\[
K_\R = \sup_{\substack{u\in H^1(\R) \\ u\not\equiv 0}} \frac{\|u\|_6^6}{\|u\|_2^4\cdot \|u'\|_2^2}
=\sup_{u\in H_\mu^1(\R) }\ \frac{\|u\|_6^6}{\mu^2\cdot \|u'\|_2^2} .
\]
Therefore
\[
\|u\|_6^6 \le K_\R \mu^2 \|u'\|_2^2 \quad \qquad \forall u\in H_\mu^1(\R).
\]
Now for every $u  \in H^1_\mu(\R)$,
\begin{align*}
6E(u,\R) &=3\|u'\|_2^2- \|u\|_6^6 \ge 3\|u'\|_2^2 - K_\R\mu^2 \|u'\|_2^2 \\
& = \|u'\|_2^2 \left(3- K_\R \mu^2\right)
\end{align*}
so that
\[
\mu^2 \le 3/K_\R\; \implies \; E(u,\R) \ge 0 \qquad\hbox{\rm for all}\; u \in H^1_\mu(\R).
\]
On the other hand, if  $\mu^2 > 3/K_\R$, and $u$ is close to optimality in the Gagliardo--Nirenberg inequality, i.e.
$\|u\|_6^6 > (K_\R-\varepsilon) \mu^2 \|u'\|_2^2$, then
\[
6E(u,\R) = 3\|u'\|_2^2 - \|u\|_6^6 \le  \|u'\|_2^2\left( 3- (K_\R-\eps) \mu^2 \right) <0
\]
for small $\varepsilon$,
and therefore
\[
\mu^2 > 3/K_\R\; \implies \; E(u,\R) <0 \qquad\hbox{\rm for some}\; u \in H^1_\mu(\R).
\]

\noindent
By mass--preserving scalings \eqref{homogeneity} it is then easy to see that
\[
\mu^2 \le 3/K_\R\; \implies \; \EE_\R(\mu) = 0,
\]
\[
\mu^2 > 3/K_\R\; \implies \; \EE_\R(\mu)  = -\infty.
\]
Therefore,
\[
\mur^2 = \frac{3}{K_\R}.
\]
This motivates the following definition.

\begin{definition}
The critical mass for a noncompact metric graph $\G$ is the number
\[
\mu_\G =  \sqrt{\frac{3}{K_\G}},
\]
where $K_\G$ is the best constant for the Gagliardo--Nirenberg inequality on $\G$.
\end{definition}

\begin{remark}
It is not difficult to see that for every noncompact $\G$,
\[
K_\R \le K_\G \le K_{\R^+}
\]
so that 
\[
\mu_{\R^+} \le \mu_\G \le \mu_\R.
\]
Thus every noncompact graph is in this sense intermediate between $\R^+$ and $\R$.
\end{remark}

\noindent
In view of the preceding discussion it is easy to prove the following statements

\begin{proposition}
\label{prop1}
Let $\G$ be a noncompact metric graph.

\begin{itemize}
\item If $\mu\le \mug$, then $\EE_\G(\mu) = 0$, and is not attained when $\mu<\mug$

\item If $\mu >\mug$, then $\EE_\G(\mu) < 0$ (possibly $-\infty$)

\item If $\mu > \mur$, then $\EE_\G(\mu) = -\infty$
\end{itemize}
\end{proposition}

\begin{corollary}
\label{coro}
A necessary condition for the existence
of a ground state of mass $\mu$ is that
\[
\mu\in[\mug,\mur].
\]
\end{corollary}

\subsection{The results}

The necessary condition of Corollary \ref{coro} is far from being sufficient.
The existence of ground states depends mainly on the topology of the graph $\G$,
according to the following  four mutually exclusive cases:
\vskip .5 cm
\begin{itemize}

\item[1.] $\quad\G$ has a terminal point (Fig. \ref{terminal})

\item[2.] $\quad\G$ satisfies Assumption (H) introduced in Sec. 4 (Fig. \ref{akka}, \ref{akkabis}).

\item[3.] $\quad\G$ has exactly one half-line and no terminal point (Fig. \ref{one})

\item[4.] $\quad\G$ has none of the above properties (Fig. \ref{none}).

\end{itemize}


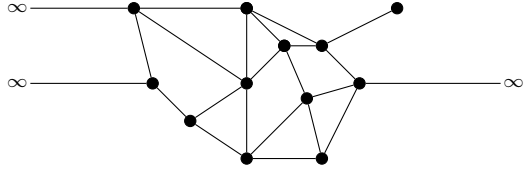
\begin{figure}
\begin{center}
\begin{tikzpicture}[xscale= 0.5,yscale=0.5]
\node at (-.5,2) [nodo] (02) {};
\node at (2,2) [nodo] (22) {};
\node at (2,4) [nodo] (24) {};
\node at (3.6,1.6) [nodo] (42) {};
\node at (3,3) [nodo] (33) {};
\node at (5,2) [nodo] (52) {};
\node at (3,3) [nodo] (32) {};
\node at (4,3) [nodo] (43) {};
\node at (-1,4) [nodo] (04) {};
\node at (2,4) [nodo] (24) {};
\node at (.5,1) [nodo] (11) {};
\node at (2,0) [nodo] (20) {};
\node at (4,0) [nodo] (40) {};
\node at (-4,2) [minimum size=0pt] (meno) {};
\node at (-4,4) [minimum size=0pt] (menoalt) {};
\node at (9,2) [minimum size=0pt] (piu) {};
\node at (-4.1,2) [infinito]  (infmeno) {$\scriptstyle\infty$};
\node at (-4.1,4) [infinito]  (infmenoalt) {$\scriptstyle\infty$};
\node at (9.1,2) [infinito]  (infpiu) {$\scriptstyle\infty$};

\node at (6,4) [nodo] (term){};
\draw [-] (43)--(term);

\draw[-] (02)--(04);
\draw[-] (04)--(24);
\draw[-] (04)--(22);
\draw[-] (24)--(22);
\draw[-] (02)--(11);
\draw[-] (11)--(22);
\draw[-] (11)--(20);
\draw[-] (20)--(22);
\draw[-] (22)--(33);
\draw[-] (24)--(33);
\draw[-] (24)--(43);
\draw[-] (33)--(43);
\draw[-] (43)--(52);
\draw[-] (33)--(42);
\draw[-] (20)--(42);
\draw[-] (20)--(40);
\draw[-] (40)--(42);
\draw[-] (42)--(52);
\draw[-] (40)--(52);
\draw[-] (02)--(meno);
\draw[-] (52)--(piu);
\draw[-] (04)--(menoalt);


\end{tikzpicture}
\end{center}
\caption{Case 1: a graph with a terminal point}
\label{terminal}
\end{figure}


%
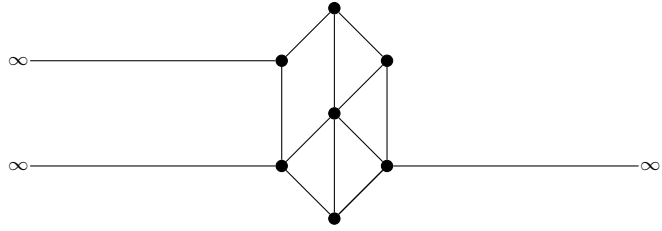
\begin{figure}[ht]
\begin{center}
\begin{tikzpicture}[xscale= 0.7,yscale=0.7]

\node at (0,0) [nodo] (00) {};
\node at (1,1) [nodo] (11) {};
\node at (0,2) [nodo] (02) {};
\node at (0,4) [nodo] (04) {};
\node at (1,3) [nodo] (13) {};
\node at (-1,1) [nodo] (-11) {};
\node at (-1,3) [nodo] (-13) {};

\node at (-6,1) [infinito] (sxbas) {$\scriptstyle\infty$};
\node at (-6,3) [infinito] (sxalt) {$\scriptstyle\infty$};
\node at (6,1) [infinito] (dx) {$\scriptstyle\infty$};

\draw[-] (00)--(11);
\draw[-] (00)--(-11);
\draw[-] (00)--(02);
\draw[-] (-11)--(02);
\draw[-] (02)--(11);
\draw[-] (00)--(11);
\draw[-] (04)--(-13);
\draw[-] (04)--(13);
\draw[-] (-13)--(-11);
\draw[-] (04)--(02);
\draw[-] (13)--(11);
\draw[-] (13)--(02);

\draw[-] (-11)--(sxbas);
\draw[-] (-13)--(sxalt);
\draw[-] (11)--(dx);


\end{tikzpicture}
\end{center}
\caption{Case 2: a graph satisfying assumption (H)}
\label{akka}
\end{figure}

\begin{figure}[ht]
\begin{center}
\begin{tikzpicture}[xscale= 0.7,yscale=0.7]

\node at (-1.4,3) [nodo] (a){};
\node at (-1.2,3) [nodo] (b){};
\node at (-1.4,1.2) [nodo] (c){};
\node at (-1.2,1.2) [nodo] (d){};
\node at (-1.4,1) [nodo] (e){};
\node at (-1.2,1) [nodo] (f){};
\node at (-.1,4) [nodo] (g){};
\node at (-.1,2.1) [nodo] (h){};
\node at (-.1,1.9) [nodo] (i){};
\node at (-.1,.1) [nodo] (l){};
\node at (-.025,-.1) [nodo] (m){};
\node at (.05,.1) [nodo] (n){};
\node at (.05,1.9) [nodo] (o){};
\node at (.2,2) [nodo] (p){};
\node at (.05,2.1) [nodo] (q){};
\node at (.05,4) [nodo] (r){};
\node at (1,3) [nodo] (s){};
\node at (1.03,2.8) [nodo] (t){};
\node at (1.03,1.2) [nodo] (u){};
\node at (1,1) [nodo] (v){};
\node at (1.2,1) [nodo] (w){};

\node at (-6,1.2) [infinito] (1) {$\scriptstyle\infty$};
\node at (-6,3) [infinito] (2) {$\scriptstyle\infty$};
\node at (-6,1) [infinito] (0) {$\scriptstyle\infty$};
\node at (6,1) [infinito] (3) {$\scriptstyle\infty$};

\draw[-,thick] (0)--(e);
\draw[-,thick] (1)--(c);
\draw[-,thick] (2)--(a);
\draw[-,thick] (3)--(w);

\draw[-,thick] (a)--(c);
\draw[-,thick] (b)--(d);
\draw[-,thick] (d)--(h);
\draw[-,thick] (h)--(g);
\draw[-,thick] (g)--(b);

\draw[-,thick] (f)--(l);
\draw[-,thick] (l)--(i);
\draw[-,thick] (i)--(f);

\draw[-,thick] (r)--(s);
\draw[-,thick] (s)--(q);
\draw[-,thick] (q)--(r);

\draw[-,thick] (p)--(u);
\draw[-,thick] (u)--(t);
\draw[-,thick] (t)--(p);

\draw[-,thick] (v)--(o);
\draw[-,thick] (o)--(n);
\draw[-,thick] (n)--(v);

\draw[-,thick] (e)--(m);
\draw[-,thick] (m)--(w);


\end{tikzpicture}
\end{center}
\caption{The graph portrayed in Fig. \ref{akka} and a covering made of $7$ cycles}
\label{akkabis}
\end{figure}
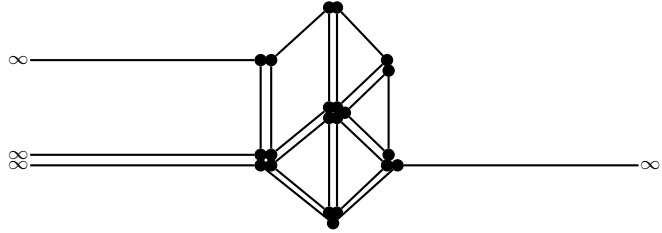

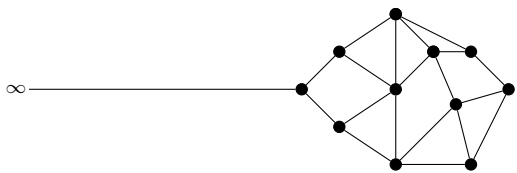
\begin{figure}[ht]
\begin{center}
\begin{tikzpicture}[xscale= 0.5,yscale=0.5]
\node at (-.5,2) [nodo] (02) {};
\node at (2,2) [nodo] (22) {};
\node at (2,4) [nodo] (24) {};
\node at (3.6,1.6) [nodo] (42) {};
\node at (3,3) [nodo] (33) {};
\node at (5,2) [nodo] (52) {};
\node at (3,3) [nodo] (32) {};
\node at (4,3) [nodo] (43) {};
\node at (.5,3) [nodo] (04) {};
\node at (2,4) [nodo] (24) {};
\node at (.5,1) [nodo] (11) {};
\node at (2,0) [nodo] (20) {};
\node at (4,0) [nodo] (40) {};
\node at (-8,2) [minimum size=0pt] (meno) {};
\node at (-8.1,2) [infinito]  (infmeno) {$\scriptstyle\infty$};

\draw[-] (02)--(04);
\draw[-] (04)--(24);
\draw[-] (04)--(22);
\draw[-] (24)--(22);
\draw[-] (02)--(11);
\draw[-] (11)--(22);
\draw[-] (11)--(20);
\draw[-] (20)--(22);
\draw[-] (22)--(33);
\draw[-] (24)--(33);
\draw[-] (24)--(43);
\draw[-] (33)--(43);
\draw[-] (43)--(52);
\draw[-] (33)--(42);
\draw[-] (20)--(42);
\draw[-] (20)--(40);
\draw[-] (40)--(42);
\draw[-] (42)--(52);
\draw[-] (40)--(52);
\draw[-] (02)--(meno);

\end{tikzpicture}
\end{center}
\caption{Case 3: a graph with exactly one half-line.}
\label{one}
\end{figure}

\hbox{}

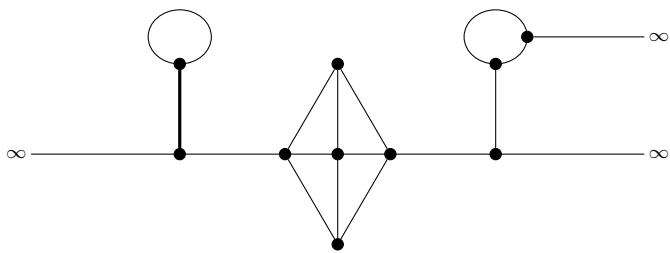
\begin{figure}[ht]
\begin{center}
\begin{tikzpicture}[xscale= 0.7,yscale=0.6]

\node at (-6,0)  [minimum size=0pt] (-60) {};
\node at (-3,0) [nodo] (-30) {};
\node at (-1,0) [nodo] (-10) {};
\node at (0,0) [nodo] (00) {};
\node at (1,0) [nodo] (10) {};
\node at (3,0) [nodo] (30) {};
\node at (6,0) [minimum size=0pt] (60) {};

\node at (0,-2) [nodo] (-02) {};
\node at (-3,2) [nodo] (-32) {};
\node at (0,2) [nodo] (02) {};
\node at (3,2) [nodo] (32) {};

\draw[-] (-60)--(60);
\draw[-] (-30)--(-32);
\draw[-,very thick] (-30)--(-32);
\draw[-] (30)--(32);
\draw[-] (-10)--(02);
\draw[-] (02)--(10);
\draw[-] (-10)--(-02);
\draw[-] (-02)--(10);
\draw[-] (-02)--(02);

\draw(-3,2.6) circle (0.6);
\draw(3,2.6) circle (0.6);
\node at (3.6,2.6) [nodo] (rig) {};
\node at (6,2.6)  [minimum size=0pt] (inf) {};
\draw[-] (rig)--(inf);

\node at (-6.1,0) [infinito] {$\scriptstyle\infty$};
\node at (6.1,0) [infinito] {$\scriptstyle\infty$};
\node at (6.1,2.6) [infinito] {$\scriptstyle\infty$};

\end{tikzpicture}
\end{center}
\caption{Case 4: a graph without terminal points, without cycle covering and with more than one half-line. No cycle can cover the thick edge.}
\label{none}
\end{figure}

\noindent
We now list the main results, case by case.

\begin{theorem}[Case 1] Assume that $\G$ has at least one terminal point. Then

\begin{itemize}
\item $\mug = \mup$
\item when $\mu \in (\mup,\mur]$, $\;\EE_\G(\mu) = -\infty$
\item when $\mu = \mup$, $\; \EE_\G(\mu) = 0$ but is attained if and only if $\G$ is a half-line.
\end{itemize}
\end{theorem}

The result shows that in the presence of a terminal point ground states do not exist (except when $\G=\R^+$).
\medskip

The terminal edge behaves like $\R^+$, almost supporting a half-soliton. The ``almost'' however cannot be eliminated,
resulting in nonexistence of ground states.

It is easy to check the second statement in the theorem. Indeed, take $u$ compactly supported on $\R^+$, with $\|u\|_2^2>\mup$,
and such that $E(u,\R^+) < 0$. The last condition can be fulfilled because the mass of $u$ is strictly larger than $\mup$, the critical mass for the half-line.
Now scale $u$ by introducing $u_\lambda(x) = \sqrt\lambda u(\lambda x)$, with $\lambda$ so large that the support of $u_\lambda$
is contained in an interval shorter than the terminal edge of $\G$.
Place $u_\lambda$ on the terminal edge of $\G$ and extend it to zero elsewhere on $\G$. Then
\[
E(u_\lambda,\G) = E(u_\lambda,\R^+)< 0,
\]
so that
\[\lim_{\lambda\to \infty}E(u_\lambda,\G)= \lim_{\lambda\to \infty}E(u_\lambda,\R^+)
= \lim_{\lambda\to \infty}\lambda^2 E(u,\R^+)= -\infty.
\]

\begin{theorem}[Case 2] Assume that $\G$ satisfies Assumption (H), so that it has a cycle covering. Then

\begin{itemize}
\item $\mug = \mur$

\item $\; \EE_\G(\mur) = 0$ and is attained if and only if $\G$ is $\R$ or a tower of bubbles.
\end{itemize}
\end{theorem}

This result shows that in the presence of a cycle covering ground states do not exist, except when $\G$ is $\R$ or a tower of bubbles (Fig. \ref{bubbles}).
\bigskip







\begin{theorem}[Case 3] 
\label{tre3}
Assume that $\G$ has exactly one half-line and no terminal point. Then
\label{tre}
\begin{itemize}
\item $\;\mug = \mup$

\item $\; \EE_\G(\mu) < 0$ (and finite) for every $\;\mu \in (\mup,\mur]$

\item $\; \EE_\G(\mu)\;$ is attained if and only if $\;\mu \in (\mup,\mur]$
\end{itemize}
\end{theorem}

This result unveils totally new phenomena: first of all,
ground states exist for a {\em whole interval} of masses, a feature
that is completely absent on the standard domains $\R$ and $\R^+$. Secondly,
ground states have {\em negative energy}, which is normal for subcritical problems,
but highly unexpected in the $L^2$--critical case.
The ultimate reason for this is the nontrivial topology of certain graphs with respect to that of $\R$ or $\R^+$.

The proof of Theorem \ref{tre}, is {\em very} involved. A sketch of some key steps will be given in Section \ref{sketch}.

We conclude with the last result, whose structure is a bit different from that of the preceding Theorems.

\begin{theorem}[Case 4] 
\label{quattro4}
Assume that $\G$ has no terminal point, no cycle covering and more than one half-line.
If, in addition,
\[
\mug < \mur,
\]
then

\begin{itemize}

\item $\; \EE_\G(\mu) < 0$ (and finite) for every $\;\mu \in (\mug,\mur]$

\item $\; \EE_\G(\mu)\;$ is attained if and only if $\;\mu \in [\mug,\mur]$
\end{itemize}
\end{theorem}

The same comments of Theorem \ref{tre3} apply:
again ground states exist for a whole interval of masses, and 
again ground states have negative energy. This time however a new feature appears: 
ground states exist also for $\mu = \mug$.
This fact is particularly interesting from the functional analytic point of view.
Indeed, since $\EE_\G(\mug) = 0$, any sequence such that $\|u_n'\|_{L^2 (\G)} \to 0$ is a minimizing sequence
and clearly compactness is lost at this level: there exist minimizing sequences at level zero that
are {\em not} precompact. However, a minimizer exists.
To obtain a ground state it is therefore necessary to select accurately a particular minimizing sequence,
in order to avoid falling onto a bad sequence.

Finally, some comments on the assumption $\mug < \mur$ are in order.

There are graphs where it is automatically satisfied, for example the {\em signpost} graph (\ref{palina}),
independently of the lengths $\ell_1, \ell_2$:
\bigskip

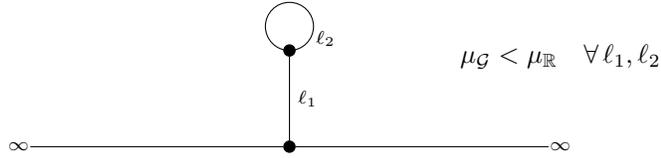
\begin{figure}
\begin{center}
\begin{tikzpicture}[scale= 0.8]

\node at (1,0) [infinito]  (5) {$\scriptstyle\infty$}; 
\node at (5.5,0) [nodo] (6) {};
\node at (10,0) [infinito]  (7) {$\scriptstyle\infty$};
\node at (5.5,1.6) [nodo] (8) {};
\draw [-] (5) -- (6) ;
\draw [-] (6) -- (7) ;
\draw [-] (6) -- (8) ;
\draw(5.5,2) circle (0.4);
\node at (5.8, .8) {$\scriptstyle\ell_1$};
\node at (6.1, 1.8) {$\scriptstyle\ell_2$};

\node at (10,1.5)  {$\mug<\mur \quad \forall\, \ell_1,\ell_2$};


\end{tikzpicture}
\end{center}
\caption{A signpost graph: for every value of $\ell_1, \ell_2$, this graph satisfies the 
hypotheses of Thm. 7.4}
\label{palina}
\end{figure}

The existence of graphs of the fourth kind where $\mug = \mur$ is an open problem.
We conjecture that in this case the sole topology of the graph is not enough to guarantee the existence of ground states.
Most likely the metric properties of the graph play a role too in this case.

\subsection{Some key steps of the existence proofs} \label{sketch}

The main ingredient of the existence proofs in Theorems \ref{tre3} and \ref{quattro4} is the following
modified Gagliardo-Nirenberg inequality, whose proof is technically very involved.

\begin{lemma}[Modified Gagliardo--Nirenberg inequality] Assume that $\G$ has no terminal point and let $\mu \le \mur$. For every $u\in H_\mu^1(\G)$ there exists $\theta\in [0,\mu]$ such that 
\begin{equation}
\label{mod}
\|u\|_{L^6(\G)}^6 \le K_\R(\mu-\theta)^2 \|u'\|_{L^2(\G)}^2 + C\theta^{1/2},
\end{equation}
with $C$ depending only on $\G$.\
\end{lemma}

We recall that the standard Gagliardo-Nirenberg inequality (with best constant) reads
\[
\|u\|_{L^6(\G)}^6 \le K_\G\mu^2 \|u'\|_{L^2(\G)}^2\qquad \forall u \in H^1_\mu(\G).
\]
In \eqref{mod} the constant $\theta$ depends on $u$. However, the inequality holds with 
a {\em smaller} constant ($K_\R \le \K_\G$), 
a {\em smaller} mass ($\mu-\theta \le \mu$) and the price is reasonable:
$C\theta^{1/2} \le C\mu^{1/2}$.

With this inequality, it is simple to show that minimizing sequences are bounded, which is the first
(and in this case more delicate) step towards an existence result.

Indeed, take $\mu \in (\mug, \mur]$. Then, by Proposition \ref{prop1}, $\EE_\G(\mu) <0$, say $\EE_\G(\mu) <-\alpha < 0$.

Let now $u_n$ be a minimizing sequence for $E$ and let $\theta_n$ be the constant in \eqref{mod} associated to $u_n$.
Then, for $n$ large,

\begin{align*}
-6\alpha \ge 6 E(u_n,\G) &= 3 \|u_n'\|_2^2 - \|u_n\|_6^6   \ge (\hbox{by \eqref{mod}}) \\ 
& \ge 3\|u_n'\|_2^2 - K_\R(\mu-\theta_n)^2\;  \|u_n'\|_2^2-C\theta_n^{1/2}\\
& \ge  \|u_n'\|_2^2\;\; (3 - K_\R(\mu -\theta_n)^2) -C\theta_n^{1/2}\\
& \ge -C\theta_n^{1/2},
\end{align*}
since 
\[
3 - K_\R(\mu -\theta_n)^2 \ge 3- K_\R\mu^2 \ge 3- K_\R\mur^2 = 0.
\]
This shows that 
\[
C\theta_n^{1/2} \ge 6\alpha,
\]
uniformly in $n$.

For this reason, $ 3 - K_\R(\mu -\theta_n)^2 \ge \delta $, so that 
\[
-6\alpha \ge \frac\delta{6}  \|u_n'\|_2^2 -C\theta_n^{1/2}
\]
from which we see that $ \|u_n'\|_2$ is uniformly bounded. Once this is established,
it is also very easy to see that $\|u_n\|_{L^p(\G)}$ is uniformly bounded too. Then one can
extract suitable subsequences and the proof of existence follows easily (in most cases).
Only in Theorem \ref{quattro4} (when $\mu=\mug$) a supplementary analysis is needed in order to 
construct a particular minimizing sequence.

\section{Conclusions and perspectives}
Inspired by physical applications, the problem of the existence of ground states for the focusing
NLS on branched structures has proved challenging from the mathematical point of view too, giving rise to a
new chapter in the Calculus of Variations, in which established techniques mix with graph theoretical notions and results. 
In particular, topology and metric of a graph interact in a highly nontrivial way, so that investigating the existence of a ground state
can involve either topological consideration or hard estimates.

The results we presented here focus on graphs with a finite number of edges and vertices, that include at least one
halfline. This means that, on the large spatial scale, all these graphs look as star graphs, and the compact core plays the role of a vertex,
possibly with an internal structure. This large-scale point of view has never been seriously considered, but it could be effective in order to 
describe ground states at low masses.

However, in all these examples the large-scale structure is still the same of a network, while it would be interesting to 
consider examples in which such a structure becomes genuinely two-dimensional, reconstructing for instance a stripe in the plane or even the entire plane.
For the first example, one could consider the case of a graph made of two parallel halflines
joined together through infinitely many parallel edges, in such a way that the distance between two consecutive edges is constant (infinite ladder graph); for the second case, the most immediate example is surely given by the 
square grid. We are currently investigating this case, and we found that the two-dimensional large-scale structure plays a very important role resulting in a
substantial change in the kind of results we are proving. Furthermore, periodicity
avoids the presence of halflines and then of quasi-solitons, so that lack of compactness in minimizing sequences 
can be  due either by spreading or by concentration only. 

\noindent{}
Beyond the problem of ground states, the issue of the existence and the shape of generic standing waves is very topical and will be addressed in forthcoming papers too.

{Far beyond these investigations, one could also think of the possibility of approximating regular domains in more dimensions with metric graphs
becoming more and more dense. This research line is very likely for a future long-term project.}


\bigskip
\begin{center}

\end{center}

\bigskip
\bigskip
\begin{minipage}[t]{10cm}
\begin{flushleft}
\small{
\textsc{Riccardo Adami}
\\*Dipartimento di Scienze Matematiche,
\\*Politecnico di Torino
\\*Corso Duca degli Abruzzi 24
\\* Torino, 10129, Italy
\\*e-mail: riccardo.adami@polito.it
\\[0.4cm]
\textsc{Enrico Serra}
\\*Dipartimento di Scienze Matematiche,
\\*Politecnico di Torino
\\*Corso Duca degli Abruzzi 24
\\* Torino, 10129, Italy
\\*e-mail: enrico.serra@polito.it
\\[0.4cm]
\textsc{Paolo Tilli}
\\*Dipartimento di Scienze Matematiche,
\\*Politecnico di Torino
\\*Corso Duca degli Abruzzi 24
\\* Torino, 10129, Italy
\\*e-mail: paolo.tilli@polito.it
}
\end{flushleft}
\end{minipage}

\end{document}